%
%
%
\newif\ifsect\newif\iffinal
\secttrue\finalfalse
\def\thm #1: #2{\medbreak\noindent{\bf #1:}\if(#2\thmp\else\thmn#2\fi}
\def\thmp #1) { (#1)\thmn{}}
\def\thmn#1#2\par{\enspace{\sl #1#2}\par
        \ifdim\lastskip<\medskipamount \removelastskip\penalty 55\medskip\fi}
\def\square{{\msam\char"03}}
\def\qedn{\thinspace\null\nobreak\hfill\square\par\medbreak}
\def\pf{\ifdim\lastskip<\smallskipamount \removelastskip\smallskip\fi
        \noindent{\sl Proof\/}:\enspace}
\def\itm#1{\item{\rm #1}\ignorespaces}

%
%
%
%
\newcount\parano
\newcount\eqnumbo
\newcount\thmno
\newcount\versiono
\newcount\remno
\newbox\notaautore
\def\neweqt#1$${\xdef #1{(\number\parano.\number\eqnumbo)}
    \eqno #1$$
    \global \advance \eqnumbo by 1}
\def\newrem #1\par{\global \advance \remno by 1
    \medbreak
{\bf Remark \the\parano.\the\remno:}\enspace #1\par
\ifdim\lastskip<\medskipamount \removelastskip\penalty 55\medskip\fi}
\def\newthmt#1 #2: #3{\xdef #2{\number\parano.\number\thmno}
    \global \advance \thmno by 1
    \medbreak\noindent
    {\bf #1 #2:}\if(#3\thmp\else\thmn#3\fi}
\def\neweqf#1$${\xdef #1{(\number\eqnumbo)}
    \eqno #1$$
    \global \advance \eqnumbo by 1}
\def\newthmf#1 #2: #3{\xdef #2{\number\thmno}
    \global \advance \thmno by 1
    \medbreak\noindent
    {\bf #1 #2:}\if(#3\thmp\else\thmn#3\fi}
\def\forclose#1{\hfil\llap{$#1$}\hfilneg}
\def\newforclose#1{
	\ifsect\xdef #1{(\number\parano.\number\eqnumbo)}\else
	\xdef #1{(\number\eqnumbo)}\fi
	\hfil\llap{$#1$}\hfilneg
	\global \advance \eqnumbo by 1
	\iffinal\else\rsimb#1\fi}
\def\forevery#1#2$${\displaylines{\let\eqno=\forclose
        \let\neweq=\newforclose\hfilneg\rlap{$\qquad\quad\forall#1$}\hfil#2\cr}$$}
\def\noNota #1\par{}
\def\today{\ifcase\month\or
   January\or February\or March\or April\or May\or June\or July\or August\or
   September\or October\or November\or December\fi
   \space\number\year}
\def\inizia{\ifsect\let\neweq=\neweqt\else\let\neweq=\neweqf\fi
\ifsect\let\newthm=\newthmt\else\let\newthm=\newthmf\fi}
\def\bititolo{\empty}
\gdef\begin #1 #2\par{\xdef\titolo{#2}
\ifsect\let\neweq=\neweqt\else\let\neweq=\neweqf\fi
\ifsect\let\newthm=\newthmt\else\let\newthm=\newthmf\fi
\iffinal\let\Nota=\noNota\fi
\centerline{\titlefont\titolo}
\if\bititolo\empty\else\medskip\centerline{\titlefont\bititolo}
\xdef\titolo{\titolo\ \bititolo}\fi
\bigskip
\centerline{\bigfont
\autore \ifvoid\notaautore\else\footnote{${}^1$}{\unhbox\notaautore}\fi}
\bigskip\if\istituto!\centerline{\today}\else
\centerline{\istituto}
\centerline{\indirizzo}
\medskip
\centerline{#1\ \anno}\fi
\bigskip\bigskip
\ifsect\else\global\thmno=1\global\eqnumbo=1\fi}
\def\anno{2011}
\def\raggedleft{\leftskip2cm plus1fill \spaceskip.3333em \xspaceskip.5em
\parindent=0pt\relax}

\font\titlefont=cmssbx10 scaled \magstep1
\font\bigfont=cmr12
\font\eightrm=cmr8
\font\sc=cmcsc10
\font\bbr=msbm10
\font\sbbr=msbm7
\font\ssbbr=msbm5
\font\msam=msam10

\font\bfm=cmmib10

\def\ca #1{{\cal #1}}

\nopagenumbers
\binoppenalty=10000
\relpenalty=10000
\newfam\amsfam
\textfont\amsfam=\bbr \scriptfont\amsfam=\sbbr \scriptscriptfont\amsfam=\ssbbr
\newfam\boldifam
\textfont\boldifam=\bfm
\let\de=\partial
\let\eps=\varepsilon
\let\phe=\varphi

\def\Ker{\mathop{\rm Ker}\nolimits}

\def\Im{\mathop{\rm Im}\nolimits}

\def\diag{\mathop{\rm diag}\nolimits}
\def\Jac{\mathop{\rm Jac}\nolimits}
\def\Span#1{\mathop{\rm Span}\nolimits\left(#1\right)}
\def\bigoperp{\mathbin{\hbox{$\bigcirc\kern-11.8pt\perp$}}}

\mathchardef\void="083F
\mathchardef\ellb="0960
\mathchardef\taub="091C
\def\C{{\mathchoice{\hbox{\bbr C}}{\hbox{\bbr C}}{\hbox{\sbbr C}}
{\hbox{\sbbr C}}}}

\def\N{{\mathchoice{\hbox{\bbr N}}{\hbox{\bbr N}}{\hbox{\sbbr N}}
{\hbox{\sbbr N}}}}

\def\Q{{\mathchoice{\hbox{\bbr Q}}{\hbox{\bbr Q}}{\hbox{\sbbr Q}}
{\hbox{\sbbr Q}}}}

\newcount\notitle
\notitle=1
\headline={\ifodd\pageno\rhead\else\lhead\fi}
\def\rhead{\ifnum\pageno=\notitle\iffinal\hfill\else\hfill\tt Version
\the\versiono; \the\day/\the\month/\the\year\fi\else\hfill\eightrm\titolo\hfill
\folio\fi}
\def\lhead{\ifnum\pageno=\notitle\hfill\else\eightrm\folio\hfill\autore\hfill
\fi}
\newbox\bibliobox
\def\setref #1{\setbox\bibliobox=\hbox{[#1]\enspace}
    \parindent=\wd\bibliobox}
\def\biblap#1{\noindent\hang\rlap{[#1]\enspace}\indent\ignorespaces}
\def\art#1 #2: #3! #4! #5 #6 #7-#8 \par{\biblap{#1}#2: {\sl #3\/}.
    #4 {\bf #5} (#6)\if.#7\else, \hbox{#7--#8}\fi.\par\smallskip}
\def\book#1 #2: #3! #4 \par{\biblap{#1}#2: {\bf #3.} #4.\par\smallskip}
\def\coll#1 #2: #3! #4! #5 \par{\biblap{#1}#2: {\sl #3\/}. In {\bf #4,}
#5.\par\smallskip}
\def\pre#1 #2: #3! #4! #5 \par{\biblap{#1}#2: {\sl #3\/}. #4, #5.\par\smallskip}
\def\raggedleft{\leftskip2cm plus1fill \spaceskip.3333em \xspaceskip.5em
\parindent=0pt\relax}
\def\Nota #1\par{\medbreak\begingroup\Bittersweet\raggedleft
#1\par\endgroup\Black
\ifdim\lastskip<\medskipamount \removelastskip\penalty 55\medskip\fi}
\newcount\defno
\def\smallsect #1. #2\par{\bigbreak\noindent{\bf #1.}\enspace{\bf #2}\par
    \global\parano=#1\global\eqnumbo=1\global\thmno=1\global\defno=0\global\remno=0
    \nobreak\smallskip\nobreak\noindent\message{#2}}
\def\newdef #1\par{\global \advance \defno by 1
    \medbreak
{\bf Definition \the\parano.\the\defno:}\enspace #1\par
\ifdim\lastskip<\medskipamount \removelastskip\penalty 55\medskip\fi}
\finalfalse
\versiono=2

\magnification\magstephalf


\let\sm=\smallskip

\finaltrue
\def\autore{Marco Abate${}^1$, Jasmin Raissy${}^2$\footnote{${}^*$}{\eightrm Partially supported by FSE, Regione Lombardia.}}
\def\indirizzo{\vbox{\hfill${}^2$Dipartimento di Matematica e Applicazioni, Universit\`a degli Studi
di Milano Bicocca,\hfill\break\null\hfill Via Cozzi 53, 20125 Milano, Italy.
E-mail: jasmin.raissy@unimib.it\hfill\null}}
\def\istituto{\vbox{\hfill${}^1$Dipartimento di Matematica, Universit\`a
di Pisa,\hfill\break\null\hfill Largo Pontecorvo 5, 56127 Pisa,
Italy. E-mail: abate@dm.unipi.it\hfill\null\vskip5pt}}

\begin {June} Formal Poincar\'e-Dulac renormalization for holomorphic germs

{{\narrower{\sc Abstract.}  Applying a general renormalization procedure for formal self-maps, producing a formal normal form simpler than the classical 
Poincar\'e-Dulac normal form, we shall give a complete 
list of normal forms for bi-dimensional superattracting germs with non-vanishing quadratic term; 
in most cases, our normal forms will be the simplest possible ones (in the sense of Wang, Zheng and Peng). We shall also discuss a few
examples of renormalization of germs tangent to the identity, revealing interesting
second-order resonance phenomena.\par
}}

\smallsect 0. Introduction

\def\autore{Marco Abate and
Jasmin Raissy}
In the study of a class of holomorphic dynamical systems, an 
important goal often is the classification under topological, holomorphic or formal 
conjugation. In particular, for each dynamical system in the class one would like to have a 
definite way of choosing a (hopefully simpler, possibly unique) representative in the same
conjugacy class; a {\sl normal form} of the original dynamical system.

The formal classification of one-dimensional germs is well-known (see, e.g., [A2]): if 
$$
f(z)=\lambda z+ a_\mu z^{\mu}+O_{\mu+1}\in\C[\![z]\!]
$$
is a one-dimensional formal power series with complex cofficients and vanishing constant term, 
where $a_\mu\ne 0$ and $O_{\mu+1}$ is a remainder term of order at least $\mu+1$, then
$f$ is formally conjugated to:
\smallskip
\item{--} $g(z)=\lambda z$ if $\lambda\ne 0$ and $\lambda$ is not a root of unity;
\item{--} $g(z)=z^\mu$ if $\lambda=0$; and to
\item{--} $g(z)=\lambda z-z^{nq+1}+\alpha z^{2nq+1}$ if $\lambda$ is a primitive $q$-th root of unity, for suitable $n\ge 1$ and $\alpha\in\C$ that are formal invariant (and $q=1$ and $n=\mu$ when $\lambda=1$).
\smallskip

In several variables, the most famous kind of normal form for local holomorphic dynamical systems
(i.e., germs of holomorphic vector fields at a singular point, or germs of 
holomorphic self-maps with a fixed point) is the {\sl Poincar\'e-Dulac normal form}
with respect to formal conjugation; let us recall very quickly its definition,
at least in the setting we are interested here, that is of formal self-maps with a fixed point,
that we can assume to be the origin in $\C^n$, without 
discussing here convergence issues.

So let $F\in\widehat\ca O^n$ be a formal transformation in $n$ complex variables,
where $\widehat\ca O^n$ denotes the space of $n$-tuples of
power series in $n$ variables  with vanishing constant term,
and let $\Lambda$ denote the (not necessarily invertible) linear term of~$F$; up to a linear change of variables, we can assume that $\Lambda$ is in Jordan normal form. For simplicity,
given a linear map $\Lambda\in M_{n,n}(\C)$ we shall denote by $\widehat\ca O^n_\Lambda$ 
the set of formal transformations in $\widehat\ca O^n$ with $\Lambda$ as
linear part.
If $\lambda_1,\ldots,\lambda_n$
are the eigenvalues of~$\Lambda$, we shall say that a multi-index $Q=(q_1,\ldots,q_n)\in\N^n$
with $q_1+\cdots+q_n\ge 2$ is {\sl $\Lambda$-resonant} if there is $j\in\{1,\ldots,n\}$
such that $\lambda_1^{q_1}\cdots\lambda_n^{q_n}=\lambda_j$. If this happens, we shall
say that the monomial $z_1^{q_1}\cdots z_n^{q_n}e_j$ is {\sl $\Lambda$-resonant,} where
$\{e_1,\ldots, e_n\}$ is the canonical basis of $\C^n$. Then (see, e.g., [Ar], [R1, 2], [R\"u])
%
given $F\in\widehat\ca O^n_\Lambda$ it is possible to find a (not unique, in general) invertible
formal transformation $\Phi\in\widehat\ca O^n_I$ with identity linear part such that
$G=\Phi^{-1}\circ F\circ\Phi$ contains only $\Lambda$-resonant monomials.

The formal transformation $G$ is a {\sl Poincar\'e-Dulac normal form} of~$F$; notice that, since
$\Phi\in\widehat\ca O^n_I$, the linear part of $G$ is still~$\Lambda$. More generally,
we shall say that a $G\in\widehat\ca O^n_\Lambda$ is {\sl in Poincar\'e-Dulac normal form}
if $G$ contains only $\Lambda$-resonant monomials.

The importance of this result cannot be underestimated, and it has been applied uncountably
many times; however it has some limitations. For instance, if $\Lambda=O$ or $\Lambda=I$
then {\it all} monomials are resonant; and thus in these cases any $F\in\widehat{\ca O}^n_\Lambda$ is in
Poincar\'e-Dulac normal form, and a further simplification (a {\sl renormalization\/}) is necessary.
Actually, even when 
a Poincar\'e-Dulac normal form is different from the original germ, it is often possible to further simplify the germ by applying invertible transformations preserving the property of being in 
Poincar\'e-Dulac normal form.  

This idea of renormalizing Poincar\'e-Dulac normal forms is not new in the context of vector 
fields; see, e.g., 
[AFGG, B1, BS, G, KOW, LS, Mu1, Mu2] and references therein.
On the other hand, with a few exceptions (see, for instance, [B2, CD]) this idea has been exploited in 
the context of self-maps only recently. One example is  [AT1], where it is applied to
a particular class of self-maps with identity linear part. More important for our
aims are [WZP1, 2], where the authors, following [KOW], construct an a priori infinite sequence
of renormalizations giving simpler and simpler normal forms.

Let us roughly describe the main ideas. For each $\nu\ge 2$ let $\ca H^\nu$ denote the space of $n$-tuples of
homogeneous polynomials in $n$ variables of degree~$\nu$. Then
every $F\in\widehat\ca O^n_\Lambda$ admits a {\sl homogeneous expansion}
$$
F=\Lambda+\sum_{\nu\ge 2} F_\nu\;,
$$
where $F_\nu\in\ca H^\nu$ is the {\sl $\nu$-homogeneous term} of~$F$. We shall also use
the notation $\{G\}_\nu$ to denote the $\nu$-homogeneous term of a formal transformation~$G$.

If $\Phi=I+\sum_{\nu\ge 2} H_\nu\in\widehat{\ca O}^n_I$
is the homogeneous expansion of an invertible formal transformation,
then it turns out that, if $L_\Lambda\colon\widehat\ca O^n\to
\widehat\ca O^n$ is defined by setting
$L_\Lambda(H)=H\circ\Lambda-\Lambda H$, then $L_\Lambda(\ca H^\nu)\subseteq\ca H^\nu$
and 
$$
\{\Phi^{-1}\circ F\circ\Phi\}_\nu=F_\nu-L_\Lambda(H_\nu)+R_\nu
\neweq\eqzPD
$$
for all $\nu\ge 2$, where $R_\nu$ is a remainder term depending only on $F_\rho$ and $H_\sigma$ with $\rho$, $\sigma<\nu$. This suggests to consider for each $\nu\ge 2$ splittings of the form
$$
\ca H^\nu=\Im L_\Lambda^\nu\oplus\ca N^\nu\quad\hbox{and}\quad
\ca H^\nu=\Ker L_\Lambda^\nu\oplus\ca M^\nu
$$ 
where $L_\Lambda^\nu= L_\Lambda|_{\ca H^\nu}$, and $\ca N^\nu$ and $\ca M^\nu$ are suitable complementary subspaces. Then \eqzPD\ implies that we can inductively choose
$H_\nu\in\ca M^\nu$ so that $\{\Phi^{-1}\circ F\circ\Phi\}_\nu\in\ca N^\nu$ for all $\nu\ge 2$;
we shall say that $G=\Phi^{-1}\circ F\circ\Phi$ is a {\sl first order normal form} of $F$ (with respect to the chosen complementary subspaces). Furthermore,
it is not difficult to see that the quadratic (actually, the first non-linear non-vanishing) homogeneous term of $G$
is a formal invariant, that is it is the same for all first order normal forms of $F$.
Notice that when $\Lambda=O$ or $\Lambda=I$ we have $L_\Lambda\equiv O$, and thus
in these cases every $F\in\widehat{\ca O}^n_\Lambda$ is a first order normal form. 

When $\Lambda$ is diagonal, $\Ker L_\Lambda$ is
generated by the resonant monomials, and $\Im L_\Lambda$ is generated by the
non-resonant monomials. Furthermore, for each $\nu\ge 2$ we have the splitting
$\ca H^\nu=\Im L_\Lambda^\nu\oplus\Ker L_\Lambda^\nu$,
and thus taking $\ca N^\mu=\Ker L_\Lambda^\nu$ and $\ca M^\mu=\Im L_\Lambda^\nu$ we have recovered the classical Poincar\'e-Dulac normal form (when $\Lambda$ has a nilpotent
part the situation is only slightly more complicated; see [Mu1, Section~4.5] for details).

Summing up, a Poincar\'e-Dulac formal normal form is composed by homogeneous terms contained in 
a complementary space of the image of the operator~$L_\Lambda$. Furthermore, 
the quadratic homogeneous term is uniquely determined, and we can still
act on the normal form by transformations having all homogeneous terms in the kernel 
of~$L_\Lambda$.

The $k$-th renormalization follows the same pattern. Assume that $F$ is in $(k-1)$-th normal form. Then there is a suitable (not necessarily
linear if $k\ge 3$; see [WZP2] for details) operator $\ca L^k$, depending on the first $k$ homogeneous terms of $F$, so that
we can bring $F$ in a normal form $G$ whose all homogeneous terms belong to a chosen complementary subspace\footnote{${}^*$}{When $k\ge 3$ one has to choose a complementary subspace to a vector space of maximal dimension contained in the image of $\ca L^k_\nu$. Actually, [WZP2] talks of ``the" subspace of maximal dimension contained in $\ca L^k_\nu$, but a priori it might not be unique.} of the image of~$\ca L^k$, and the first $k+1$ homogeneous terms
of $G$ are uniquely determined; we shall say that $G$ is in {\sl $k$-th order normal form}
(with respect to the chosen subspaces).

A formal transformation $G$ is in {\sl  infinite order normal form} if it is in $k$-th normal form 
for all $k$, with respect to some choice of complementary subspaces and using the operators 
$\ca L^k$ defined using the first $k$ homogeneous terms of~$G$. The main result of
[WZP2] then states that every element of $\widehat{\ca O}^n_\Lambda$ can be brought to
a (possibly not unique) infinite order normal form by a sequence of formal conjugations tangent to the identity.

To apply these results, we need a rule for choosing complementary subspaces. It turns out
that an efficient way of doing this is by taking
orthogonal complements with respect to the Fischer Hermitian product, defined by (see [F])
$$
\langle z_1^{p_1}\cdots z_n^{p_n}e_h,z_1^{q_1}\cdots z_n^{q_n}e_k\rangle=
\cases{
0&if $h\ne k$ or $p_j\ne q_j$ for some $j$;\cr
\noalign{\smallskip}
\displaystyle {p_1!\cdots p_n!\over(p_1+\cdots+p_n)!}&if $h=k$ and $p_j=q_j$ for all $j$.\cr}
\neweq\eqzFi
$$
With this choice, as we shall see in Sections~2 and~3,  
the expression of the second order (and often infinite order) normal forms can be quite simple. For instance, in Section~2 we shall apply this procedure to the case of superattracting
(i.e., with $\Lambda=O$) 2-dimensional formal transformations, case that has no analogue
in the vector field setting, proving the following

\newthm Theorem \zuno: Let $F\in\widehat{\ca O}^2_O$ be of the form $F(z,w)=F_2(z,w)+O_3$.
Then:
\smallskip
\itm{(i)} if $F_2(z,w)=(z^2,zw)$ or $F_2(z,w)=(-z^2,-z^2-zw)$ then $F$ is formally conjugated to an unique infinite order normal form
$$
G(z,w)=F_2(z,w)+\bigl(\phe(w)-z\psi'(w),2\psi(w)\bigr)\;,
$$
where $\phe$, $\psi\in\C[\![\zeta]\!]$ are power series of order at least~3;
\itm{(ii)} if $F_2(z,w)=(-zw,-z^2-w^2)$ then 
$F$ is formally conjugated to an unique infinite order normal form
$$
G(z,w)=F_2(z,w)+\bigl(-2\phe(z+w)+2\psi(w-z),\phe(z+w)+\psi(w-z)\bigr)\;,
$$
where $\phe$, $\psi\in\C[\![\zeta]\!]$ are power series of order at least~3;
\itm{(iii)} if $F_2(z,w)=(zw,zw+w^2)$ then 
$F$ is formally conjugated to an unique infinite order normal form
$$
G(z,w)=F_2(z,w)+\bigl(w\phe'(z)+\psi(z),2\phe(z)-w\phe'(z)-\psi(z)\bigr)\;,
$$
where $\phe$, $\psi\in\C[\![\zeta]\!]$ are power series of order at least~3;
\itm{(iv)} if $F_2(z,w)=\bigl(-\rho z^2,(1-\rho)zw\bigr)$ with $\rho\ne 0$,~$1$ then 
$F$ is formally conjugated to an unique infinite order normal form
$$
G(z,w)=F_2(z,w)+\bigl((\rho-1)z\phe'(w)+\psi(w),-2\rho\phe(z)\bigr)\;,
$$
where $\phe$, $\psi\in\C[\![\zeta]\!]$ are power series of order at least~3;
\itm{(v)} if $F_2(z,w)=(-z^2+zw,w^2\bigr)$ then 
$F$ is formally conjugated to an unique infinite order normal form
$$
G(z,w)=F_2(z,w)+\bigl(\phe({\textstyle{z\over 2}}+w),-{\textstyle{1\over 4}}\phe({\textstyle{z\over 2}}+w)+\psi(z)\bigr)\;,
$$
where $\phe$, $\psi\in\C[\![\zeta]\!]$ are power series of order at least~3;
\itm{(vi)} if $F_2(z,w)=\bigl(\rho z^2+zw,(1+\rho)zw+w^2\bigr)$ with $\rho\ne 0$,~$-1$ then 
$F$ is formally conjugated to an unique infinite order normal form
$$
\eqalign{
G(z,w)=F_2(z,w)+&\left({1\over\rho}\left[{1-\sqrt{-\rho}\over 2m_\rho^2}\phe(m_\rho z+w)
+{1+\sqrt{-\rho}\over 2n_\rho^2}\phe(n_\rho z+w)\right]\right.\cr
&\phantom{\biggl(\,}
+{1+\rho\over2\sqrt{-\rho}}\left(
{1\over m_\rho^2}\psi(m_\rho z+w)-{1\over n_\rho^2}\psi(n_\rho z+w)\right),\cr
&\quad {1-\sqrt{-\rho}\over2}\phe(m_\rho z+w)+{1+\sqrt{-\rho}\over2}\phe(n_\rho z+w)\cr
&\phantom{\biggl(\;}+\left.{\rho(1+\rho)\over2\sqrt{-\rho}}\bigl(\psi(m_\rho z+w)-\psi(n_\rho z+w)\bigr)
\right)\cr}
$$ 
where $\sqrt{-\rho}$ is any square root of $-\rho$,
$$
m_\rho={\sqrt{-\rho}-\rho\over \rho(1+\rho)}\;,\quad n_\rho=-{\sqrt{-\rho}+\rho\over\rho(1+\rho)}\;,
$$
and $\phe$, $\psi\in\C[\![\zeta]\!]$ are power series of order at least~3;
\itm{(vii)} if 
$F_2(z,w)=\bigl(\rho (-z^2+zw),(1-\rho)(zw-w^2)\bigr)$
with $\rho\ne 0$,~$1$ then $F$ is formally conjugated to an unique infinite order normal form
$$
\eqalign{
G(z,w)=F_2(z,w)+&\left(z{\de\over\de z}\bigl[\phe(z+w)+\psi(z+w)\bigr]-\phe(z+w),\right.\cr
&\left.\quad{\rho-1\over\rho}\left(z{\de\over\de z}\bigl[\phe(z+w)-\psi(z+w)\bigr]
-3\phe(z+w)+2\psi(z+w)\right)
\right)\cr}
$$ 
where $\phe$,~$\psi\in\C[\![\zeta]\!]$ are power series of order at least~3;
\itm{(ix)} if $F_2(z,w)=(-z^2,-w^2)$ then $F$ is formally conjugated to an unique infinite order normal form
$$
G(z,w)=F_2(z,w)+\bigl(\phe(w),\psi(z)\bigr)
$$
where $\phe$,~$\psi\in\C[\![\zeta]\!]$ are power series of order at least~3;
\itm{(x)} if $F_2(z,w)=(-\rho z^2,(1-\rho)zw-w^2)$ with $\rho\ne0$,~$1$ then $F$ is formally conjugated to an unique infinite order normal form
$$
G(z,w)=F_2(z,w)+\left(\phe(w)+{(1-\rho)^2\over 4\rho}\psi\left({2\over 1-\rho}z+w\right),\psi
\left({2\over1-\rho}z+w\right)\right)
$$
where $\phe$,~$\psi\in\C[\![\zeta]\!]$ are power series of order at least~3;
\itm{(xi)} if $F_2(z,w)=(-\rho z^2+(1-\tau)zw,(1-\rho)zw-\tau w^2)$ with $\rho$,~$\tau\ne0$,~$1$ 
and $\rho+\tau\ne 1$ then $F$ is formally conjugated to an unique infinite order normal form
$$
\eqalign{
G(z,w)=F_2(z,w)+&\left(
{\tau\over\rho}\left[{\sqrt{\rho+\tau-1}+\sqrt{\rho\tau}\over 2m_{\rho,\tau}^2}\phe(m_{\rho,\tau}z+w)
\right.\right.\cr
&\phantom{{\tau\over\rho}\biggl[\quad}+{\sqrt{\rho+\tau-1}-\sqrt{\rho\tau}\over 2n_{\rho,\tau}^2}\phe(n_{\rho,\tau}z+w)\cr
&\phantom{{\tau\over\rho}\biggl[\quad}\left.+
{1\over m_{\rho,\tau}^2}\psi(m_{\rho,\tau}z+w)-{1\over n_{\rho,\tau}^2}\psi(n_{\rho,\tau}z+w)\right],
\cr
&\quad{\sqrt{\rho+\tau-1}+\sqrt{\rho\tau}\over 2}\phe(m_{\rho,\tau}z+w)\cr
&\phantom{\quad\,}
+{\sqrt{\rho+\tau-1}-\sqrt{\rho\tau}\over 2}\phe(n_{\rho,\tau}z+w)\cr
&\phantom{\quad\,}+\psi(m_{\rho,\tau}z+w)-\psi(n_{\rho,\tau}z+w)
\Biggr)\;,\cr}
$$ 
where
$$
m_{\rho,\tau}={\sqrt{\rho\tau}\sqrt{\rho+\tau-1}-\rho\tau\over\rho(\rho-1)}\;,\qquad n_{\rho,\tau}
=-{\sqrt{\rho\tau}\sqrt{\rho+\tau-1}+\rho\tau\over\rho(\rho-1)}\;.
$$
\indent and 
$\phe$,~$\psi\in\C[\![\zeta]\!]$ are power series of order at least~3.

In [A1] we showed that the list of quadratic terms
in this theorem gives a complete list of all possible quadratic terms up to linear change of coordinates, with the exception of four degenerate cases where one of the coordinates is identically zero. In these cases missing we shall anyway be able to give 
a second order normal form:

\newthm Proposition \zdue: Let $F\in\widehat{\ca O}^2_O$ be of the form $F(z,w)=F_2(z,w)+O_3$.
Then:
\smallskip
\itm{(i)} if $F_2(z,w)=(0,-z^2)$ then $F$ is formally conjugated to a unique second order normal form
$$
G(z,w)=F_2(z,w)+\bigl(\Phi(z,w),\psi(w)\bigr)\;,
$$
where $\psi\in\C[\![\zeta]\!]$ and $\Phi\in\C[\![z,w]\!]$ are power series of order at least~3;
\itm{(ii)} if $F_2(z,w)=(0,zw)$ then $F$ is formally conjugated to a unique second order normal form
$$
G(z,w)=F_2(z,w)+\bigl(\Phi(z,w),0\bigr)\;,
$$
where $\Phi\in\C[\![z,w]\!]$ is a power series of order at least~3;
\itm{(iii)} if $F_2(z,w)=(-z^2,0)$ then $F$ is formally conjugated to a unique second order normal form
$$
G(z,w)=F_2(z,w)+\bigl(\psi(w),\Phi(z,w)\bigr)\;,
$$
where $\psi\in\C[\![\zeta]\!]$ and $\Phi\in\C[\![z,w]\!]$ are power series of order at least~3;
\itm{(iv)} if $F_2(z,w)=(z^2-zw,0)$ then $F$ is formally conjugated to a unique second order normal form
$$
G(z,w)=F_2(z,w)+\bigl(0,\Phi(z,w)\bigr)\;,
$$
\indent where $\Phi\in\C[\![z,w]\!]$ is a power series of order at least~3.

Finally, in Section~3 we shall also discuss a few interesting examples with $\Lambda=I$,
showing in particular the appearance of non-trivial second-order resonance phenomena. For instance, we shall prove the following

\newthm Proposition \ztre: Let $F\in\widehat{\ca O}^2_I$ be of the form $F(z,w)=(z,w)+F_2(z,w)+O_3$, with 
$$
F_2(z,w)=\bigl(-\rho z^2,(1-\rho)zw\bigr)
$$ 
and $\rho\ne 0$. Put
$$
\ca E=\bigl([0,1]\cap\Q\bigr)\cup\left\{-{1\over n}\biggm| n\in\N^*\right\}
\quad\hbox{and}\quad
\ca F=\bigl([0,1]\cap\Q\bigr)\cup\left\{1+{1\over n},1+{2\over n}\biggm| n\in\N^*\right\}\;.
$$
Then:
\smallskip
\item{\rm(i)} if $\rho\notin\ca E\cup\ca F$ then $F$ is formally conjugated to a unique second order normal form
$$
G(z,w)=(z,w)+F_2(z,w)+\bigl(az^3+\phe(w)+(1-\rho)z\psi'(w),(1-\rho)w\psi'(w)+(3\rho-1)\psi(z)\bigr)\;,
$$
where $\phe$,~$\psi\in\C[\![\zeta]\!]$ are power series of order at least~3, and $a\in\C$;
\item{\rm(ii)} if $\rho=1+{1\over n}\in\ca F\setminus\ca E$ then $F$ is formally conjugated to a unique second order normal form
$$
\eqalign{
G(z,w)=(z,w)&+F_2(z,w)\cr
&+\left(a_0z^3+a_1z^2w^{n+1}+\phe(w)-{1\over n}z\psi'(w),-{1\over n}w\psi'(w)+\left(2+{3\over n}\right)\psi(w)\right)\;,\cr}
$$
where $\phe$,~$\psi\in\C[\![\zeta]\!]$ are power series of order at least~3, and $a_0$,~$a_1\in\C$;
\item{\rm(iii)} if $\rho=1+{2\over m}\in\ca F\setminus\ca E$ with $m$ odd then $F$ is formally conjugated to a unique second order normal form
$$
\eqalign{
G(z,w)=(z,w)&+F_2(z,w)\cr
&+\left(a_0z^3+\phe(w)-{2\over m}z\bigl(w\psi'(w)+\psi(w)\bigr),\right.\cr
&\qquad\quad\left.-{2\over m}w^2\psi'(w)+\left(2+{4\over m}\right)w\psi(w)\right)\;,\cr}
$$
where $\phe$,~$\psi\in\C[\![\zeta]\!]$ are power series of order at least respectively~3 and~2, and $a_0\in\C$;
\item{\rm(iv)} if $\rho=-{1\over n}\in\ca E\setminus\ca F$ then $F$ is formally conjugated to a unique second order normal form
$$
\eqalign{
G(z,w)=(z,w)&+F_2(z,w)\cr
&+\left(a_0z^3+\phe(w)+\left(1+{1\over n}\right)z\bigl(w\psi'(w)+\psi(w)\bigr),\right.\cr
&\qquad\quad\left.a_1z^{n+2}+\psi(z)+\left(1+{1\over n}\right)w^2\psi'(w)-{2\over n}w\psi(w)\right)\;,\cr}
$$
where $\phe$,~$\psi\in\C[\![\zeta]\!]$ are power series of order at least respectively~3 and~2, and $a_0$~$a_1\in\C$;
\item{\rm(v)} if $\rho=1\in\ca E\cap\ca F$ then $F$ is formally conjugated to a unique second order normal form
$$
G(z,w)=(z,w)+F_2(z,w)+\left(\phe_1(w)+z^3\psi(w),\phe_2(w)+z\phe_3(w)\right)\;,
$$
where $\phe_1$,~$\phe_2\in\C[\![\zeta]\!]$ are power series of order at least 3,
$\phe_2\in\C[\![\zeta]\!]$ is a power series of order at least~2, and $\phe_3\in\C[\![\zeta]\!]$ is a power series;
\item{\rm(vi)} if $\rho=a/b\in(0,1)\cap\Q\subset\ca E\setminus\ca F$ then $F$ is formally conjugated to a unique second order normal form
$$
\eqalign{
G(z,w)&=(z,w)+F_2(z,w)\cr
&+\!\!\left(\phe(w)+z^3\phe_0(z^{b-a}w^a)+(b-a){\de\over\de w}\bigl(z^2w\chi(z^{b-a}w^a)\bigr)\!+\!
\left(1-{a\over b}\right)\!z\bigl(w\psi'(w)+\psi(w)\bigr),\right.\cr
&\qquad\quad\left.
a{\de\over\de z}\bigl(z^2w\chi(z^{b-a}w^a)\bigr)+\left(1-{a\over b}\right)w^2\psi'(w)+2{a\over b}
w\psi(w)\right)\;,\cr}
$$
\indent where $\phe$,~$\psi\in\C[\![\zeta]\!]$ are power series of order at least~3, and $\phe_0$,~$\chi\in\C[\![\zeta]\!]$ are power series of\break\indent order at least~1.

\smallsect 1. Renormalization

In this section we shall recover, with a different proof, the part of the renormalization procedure of [WZP2] useful for our aims. One difference between our approach and theirs is that we shall
systematically use the relations between homogeneous polynomials and symmetric multilinear
maps instead of relying on higher order derivatives as in [WZP2].

Let us start collecting a few results on homogeneous polynomials and maps we shall need later.

\newdef We shall denote by  $\ca H^d$ the 
space of {\sl homogenous maps of degree $d$,} i.e., of $n$-tuples of homogeneous polynomials of degree~$d\ge 1$ in the variables
$(z_1,\ldots,z_n)$. It is well known (see, e.g., [C, pp. 79--88]) that to each $P\in\ca H^d$ is
associated a {\it unique} symmetric multilinear map $\tilde P\colon(\C^n)^d\to\C^n$ such that
$$
P(z)=\tilde P(z,\ldots,z)
$$
for all $z\in\C^n$. We also set $\ca H=\prod\limits_{d\ge 2}\ca H^d$.


Roughly speaking, the symmetric multilinear map associated to a homogeneous map~$H$ encodes the derivatives of~$H$. For instance, it is easy to check that for each $H\in\ca H^d$ we have
$$
(\Jac H)(z)\cdot v=d\,\tilde H(v,z,\ldots,z)
\neweq\equuno
$$ 
for all $z$, $v\in\C^n$.


Later on we shall need to compute the multilinear map associated to a homogeneous map
obtained as a composition. The formula we are interested in is contained in the next lemma.

\newthm Lemma \udueb: Assume that $P\in\ca H^d$ is of the form
$$
P(z)=\tilde K\bigl(H_{d_1}(z),\ldots,H_{d_r}(z)\bigr)\;,
$$
where $\tilde K$ is $r$-multilinear, $d_1+\cdots+d_r=d$, and $H_{d_j}\in\ca H^{d_j}$ for
$j=1,\ldots, r$. Then
$$
\tilde P(v,w,\ldots,w)={1\over d}\sum_{j=1}^r d_j\tilde K\bigl(H_{d_1}(w),\ldots,\tilde H_{d_j}(v,w,\ldots,w),
\ldots,H_{d_r}(w)\bigr)
$$
for all $v$, $w\in\C^n$.

\pf Write $z=w+\eps v$. Then
$$
\eqalign{P(w)&+d\eps\tilde P(v,w,\ldots,w)+O(\eps^2)\cr
&=P(w+\eps v)=
\tilde K\bigl(\tilde H_{d_1}(w+\eps v,\ldots,w+\eps v),\ldots,\tilde H_{d_r}(w+\eps v,\ldots,
w+\eps v)\bigr)\cr
&=\tilde K\bigl(H_{d_1}(w),\ldots,H_{d_r}(w)\bigr)\!+\eps\sum_{j=1}^r \!d_j
\tilde K\bigl(H_{d_1}(w),\ldots,\tilde H_{d_j}(v,w,\ldots,w),\ldots,H_{d_r}(w)\bigr)+O(\eps^2)\;,
\cr}
$$
and we are done.\qedn

\newdef Given a linear map $\Lambda\in M_{n,n}(\C)$, we define a linear operator 
$L_\Lambda\colon\ca H\to\ca H$ by setting
$$
L_\Lambda(H)=H\circ\Lambda -\Lambda H\;.
$$
We shall say that a homogeneous map $H\in\ca H^d$ is {\sl $\Lambda$-resonant} if $L_\Lambda(H)=O$, and we shall denote by $\ca H^d_\Lambda=\Ker L_\Lambda\cap \ca H^d$ the 
subspace of $\Lambda$-resonant homogeneous maps of degree~$d$. Finally,
we set $\ca H_\Lambda=\prod\limits_{d\ge 2}\ca H^d_\Lambda$. 

When $\Lambda$ is diagonal, then the $\Lambda$-resonant monomials are exactly the
resonant monomials appearing in the classical Poincar\'e-Dulac theory.

\newdef If $Q=(q_1,\ldots,q_n)\in\N^n$ is a multi-index
and $z=(z_1,\ldots,z_n)\in\C^n$, we shall put $z^Q=z_1^{q_1}\cdots z_n^{q_n}$. Given
$\Lambda=\diag(\lambda_1,\ldots,\lambda_n)\in M_{n,n}(\C)$, we shall say that $Q\in\N^n$ with $q_1+\cdots+q_n\ge 2$
is {\sl $\Lambda$-resonant on the $j$-th coordinate} if $\lambda_1^{q_1}\cdots
\lambda_n^{q_n}=\lambda_j$. If $Q$ is $\Lambda$-resonant on the $j$-th coordinate, we shall
also say that the monomial $z^Q e_j$ is {\sl $\Lambda$-resonant,} where $\{e_1,\ldots,e_n\}$
is the canonical basis of~$\C^n$.

\newrem If $\Lambda=\diag(\lambda_1,\ldots,\lambda_n)\in M_{n,n}(\C)$ is diagonal, and $z^Qe_j
\in\ca H^d$ is a homogeneous monomial (with $q_1+\cdots+q_n=d$), then (identifying the matrix $\Lambda$ with the vector, still denoted by $\Lambda$, of its diagonal entries) we have 
$$
L_\Lambda(z^Qe_j)=(\Lambda^Q-\lambda_j)z^Qe_j\;.
$$
Therefore $z^Qe_j$ is $\Lambda$-resonant if and only if $Q$ is $\Lambda$-resonant in the $j$-th coordinate. In particular, a basis of $\ca H^d_\Lambda$ is given by the $\Lambda$-resonant
monomials, and we have
$$
\ca H^d=\ca H^d_\Lambda\oplus\Im L_\Lambda|_{\ca H^d}
$$
for all $d\ge 2$.

It is possible to detect the $\Lambda$-resonance by using the associated multilinear map:

\newthm Lemma \utre: If $\Lambda\in M_{n,n}(\C)$ and $H\in\ca H^d$ then $H$ is
$\Lambda$-resonant if and only if
$$
\tilde H(\Lambda v_1,\ldots,\Lambda v_d)=\Lambda\tilde H(v_1,\ldots,v_d)
\neweq\equdue
$$
for all $v_1,\ldots,v_d\in\C^n$. In particular, if $H\in\ca H^d_\Lambda$ then
$$
\bigl((\Jac H)\circ\Lambda\bigr)\cdot\Lambda=\Lambda\cdot(\Jac H)\;.
\neweq\equtre
$$

\pf One direction is trivial. Conversely, assume $H\in\ca H^d_\Lambda$. By definition, $H$ is $\Lambda$-resonant if and only if $\tilde H(\Lambda w,\ldots,\Lambda w)
=\Lambda\tilde H(w,\ldots,w)$ for all $w\in\C^n$. Put $w=z+\eps v_1$; then
$$
\eqalign{\tilde H(\Lambda z,\ldots,\Lambda z)+\eps d\,\tilde H(\Lambda v_1,\Lambda z,\ldots,
\Lambda z)+O(\eps^2)&=\tilde H\bigl(\Lambda(z+\eps v_1),\ldots,\Lambda(z+\eps v_1)\bigr)\cr
&=\Lambda\tilde H(z+\eps v_1,\ldots,z+\eps v_1)\cr
&=\Lambda\tilde H(z,\ldots,z)+\eps d\,\Lambda\tilde H(v_1,z,\ldots,z)+O(\eps^2)\;,\cr}
$$
and thus
$$
\tilde H(\Lambda v_1,\Lambda z,\ldots,\Lambda z)=\Lambda\tilde H(v_1,z,\ldots,z)\;;
\neweq\equdueb
$$
in particular \equtre\ is a consequence of \equuno.

Now put $z=z_1+\eps v_2$ in \equdueb. We get
$$
\eqalign{\tilde H(\Lambda v_1,\Lambda z_1,\ldots,\Lambda z_1)&+\eps (d-1)\tilde H(\Lambda v_1,\Lambda v_2,\Lambda z_1,\ldots,
\Lambda z_1)+O(\eps^2)\cr
&=\tilde H\bigl(\Lambda v_1,\Lambda(z_1+\eps v_2),\ldots,\Lambda(z_1+\eps v_2)\bigr)\cr
&=\Lambda\tilde H(v_1,z_1+\eps v_2,\ldots,z_1+\eps v_2)\cr
&=\Lambda\tilde H(v_1,z_1,\ldots,z_1)+\eps (d-1)\Lambda\tilde H(v_1,v_2,z,\ldots,z)+O(\eps^2)\;,\cr}
$$
and hence
$$
\tilde H(\Lambda v_1,\Lambda v_2,\Lambda z_1,\ldots,\Lambda z_1)=\Lambda\tilde H(v_1,v_2,z_1,\ldots,z_1)
$$
for all $v_1$, $v_2$, $z_1\in\C^n$. Proceeding in this way we get \equdue.
\qedn

We now introduce the operator needed for the second order normalization.

\newdef Given $P\in\ca H^\mu$ and $\Lambda\in M_{n,n}(\C)$, let $L_{P,\Lambda}\colon
\ca H^d\to\ca H^{d+\mu-1}$ be given by
$$
L_{P,\Lambda}(H)(z)=d\,\tilde H\bigl(P(z),\Lambda z,\ldots,
\Lambda z\bigr)-\mu\tilde P\bigl(H(z),z,\ldots,z\bigr)\;.
$$

\newrem Equation \equuno\ implies that
$$
d\,\tilde H\bigl(P(z),\Lambda z,\ldots,\Lambda z\bigr)=(\Jac H)(\Lambda z)\cdot P(z)\;.
$$
Therefore
$$
L_{P,\Lambda}(H)=\bigl((\Jac H)\circ\Lambda\bigr)\cdot P-(\Jac P)\cdot H\;;
$$
In the notations of [WZP2] we have $L_{P,\Lambda}(H)=[H,P)$, and $L_{P,\Lambda}|_{\ca H^d_\Lambda}=\ca T_d[P]$ when $P\in\ca H^\mu_\Lambda$.


Using multilinear maps it is easy to prove the following useful fact (cp. [WZP2, Lemma~2.1]):

\newthm Lemma \uqua: Take $\Lambda\in M_{n,n}(\C)$ and $P\in\ca H^\mu_\Lambda$.
Then $L_{P,\Lambda}(\ca H^d_\Lambda)\subseteq\ca H^{d+\mu-1}_\Lambda$ for all $d\ge 2$.

\pf Using Lemma~\utre\ and the definition of $L_{P,\Lambda}$, if $H\in\ca H^d_\Lambda$ we get
$$
\eqalign{
L_{P,\Lambda}(H)(\Lambda z)&=d\,\tilde H\bigl(P(\Lambda z),\Lambda^2 z,\ldots,\Lambda^2 z\bigr)
-\mu\tilde P\bigl(H(\Lambda z),\Lambda z,\ldots,\Lambda z\bigr)\cr
&=d\,\tilde H\bigl(\Lambda P(z),\Lambda^2 z,\ldots,\Lambda^2 z\bigr)
-\mu\tilde P\bigl(\Lambda H(z),\Lambda z,\ldots,\Lambda z\bigr)\cr
&=d\,\Lambda\tilde H\bigl(P(z),\Lambda z,\ldots,\Lambda z\bigr)
-\mu\Lambda \tilde P\bigl(H(z),z,\ldots,z\bigr)\cr
&=\Lambda L_{P,\Lambda}(H)(z)\;.
\cr}
$$ 
\qedn

%
%
%

To state and prove the main technical result of this section we fix a few more notations.

\newdef We shall denote by $\widehat\ca O^n=\prod\limits_{d\ge 1}\ca H^d$
the space of $n$-tuples 
of formal power series with vanishing constant term. Furthermore, given $\Lambda\in M_{n,n}(\C)$ we shall denote by $\widehat\ca O^n_\Lambda$ the subset of~$F\in\widehat\ca O^n$ with $dF_O=\Lambda$. Every $F\in\widehat\ca O^n$ can be written in a unique way
as a formal sum
$$
F=\sum_{d\ge 1} F_d
\neweq\eqduno
$$
with $F_d\in\ca H^d$; \eqduno\ is the {\sl homogeneous expansion} of $F$, and $F_d$ 
is the {\sl $d$-homogeneous term} of $F$. We shall often write $\{F\}_d$ for $F_d$. In particular,
if $F\in\widehat\ca O^n_\Lambda$ then $\{F\}_1=\Lambda$.

The homogeneous terms behave in a predictable way with respect to composition and inverse:
indeed it is easy to see that if
$F=\sum\limits_{d\ge 1} F_d$ and
$G=\sum\limits_{d\ge 1} G_d$ are two elements of $\widehat{\ca O}^n$ then
$$
\{F\circ G\}_d=\sum_{{1\le r\le d}\atop{d_1+\cdots+d_r=d}}\tilde F_r(G_{d_1},\ldots, G_{d_r})
\neweq\eqdcomp
$$
for all $d\ge 1$; and that if $\Phi=I+\sum\limits_{d\ge2}H_d$ belongs to $\widehat{\ca O}^n_I$
then the homogeneous expansion of the inverse transformation $\Phi^{-1}=I+\sum\limits_{d\ge 2}K_d$ is given by
$$
K_d=
-H_d-\sum_{{2\le r\le d-1}\atop{d_1+\cdots+d_r=d}}\tilde K_r(H_{d_1},\ldots,H_{d_r})
\neweq\eqddue
$$
for all $d\ge 2$. In particular we have

\newthm Lemma \ddue: Let $\Phi=I+\sum\limits_{d\ge2}H_d\in\widehat\ca O^n_I$, and let $\Phi^{-1}=I+\sum\limits_{d\ge 2}K_d$
be the homogeneous expansion of the inverse. Then if $H_2,\ldots, H_d$ are $\Lambda$-resonant for some $\Lambda\in M_{n,n}(\C)$ and $d\ge 2$
then also $K_d$ is. 

\pf We argue by induction. Assume that $H_2,\ldots, H_d$ are
$\Lambda$-resonant. If $d=2$ then $K_2=-H_2$ and
thus $K_2$ is clearly $\Lambda$-resonant. Assume the assertion true for $d-1$; in particular,
$K_2,\ldots,K_{d-1}$ are $\Lambda$-resonant. Then
$$
\eqalign{K_d\circ\Lambda&=-H_d\circ\Lambda-\sum_{{2\le r\le d-1}\atop{d_1+\cdots+d_r=d}}\tilde K_r(H_{d_1}\circ\Lambda,\ldots,H_{d_r}\circ\Lambda)\cr
&=\Lambda H_d-\sum_{{2\le r\le d-1}\atop{d_1+\cdots+d_r=d}}\tilde K_r(\Lambda H_{d_1},\ldots,\Lambda H_{d_r})=\Lambda K_d\cr}
$$
because $K_2,\ldots,K_{d-1}$ are $\Lambda$-resonant (and we are using Lemma~\utre).
\qedn

\newdef Given $\Lambda\in M_{n,n}(\C)$, we shall say that $F\in\widehat\ca O^n$ is 
{\sl $\Lambda$-resonant} if $F\circ\Lambda=\Lambda F$. Clearly, $F$ is $\Lambda$-resonant
if and only if $\{F\}_d\in\ca H^d_\Lambda$ for all $d\in\N$. 

%
%
%

The main technical result of this section is the following analogue of [WZP2, Theorem~2.4]:

\newthm Theorem \dquattro: Given $F\in\widehat\ca O^n_O$, let
$F=\Lambda+\sum\limits_{d\ge\mu}F_d$ be its homogeneous expansion, with $F_\mu\ne O$.
Then for every $\Phi\in\widehat\ca O^n_I$ with homogeneous
expansion $\Phi=I+\sum\limits_{d\ge 2}H_d$ and every $\nu\ge 2$ we have
$$
\{\Phi^{-1}\circ F\circ\Phi\}_\nu=F_\nu-L_\Lambda(H_\nu)-L_{F_\mu,\Lambda}(H_{\nu-\mu+1})+Q_\nu+R_\nu\;,
\neweq\eqdqua
$$
where $Q_\nu$ depends only on $\Lambda$ and on $H_\gamma$ with $\gamma<\nu$, while $R_\nu$ depends only on $F_\rho$ with $\rho<\nu$ and on $H_\gamma$ with $\gamma<\nu-\mu+1$, and we put $L_{F_\mu,\Lambda}(H_1)=O$. Furthermore,
we have:
\smallskip
\itm{(i)} if $H_2,\ldots,H_{\nu-1}\in\ca H_\Lambda$ then $Q_\nu=O$; in particular, if $\Phi$ is $\Lambda$-resonant then $L_\Lambda(H_\nu)=Q_\nu=O$
for all $\nu\ge 2$;
\itm{(ii)} if $\Phi$ is $\Lambda$-resonant then $\{\Phi^{-1}\circ F\circ\Phi\}_\nu=O$ for $2\le \nu<\mu$, $\{\Phi^{-1}\circ F\circ\Phi\}_\mu=
F_\mu$, and 
$$
\{\Phi^{-1}\circ F\circ\Phi\}_{\mu+1}=F_{\mu+1}-L_{F_\mu,\Lambda}(H_2)\;;
$$
\itm{(iii)} if $F=\Lambda$ then $R_\nu=O$ for all $\nu\ge 2$;
\item{\rm(iv)} if $F_2,\ldots, F_{\nu-1}$ and $H_2,\ldots,H_{\nu-\mu}$ are $\Lambda$-resonant
then $R_\nu$ is $\Lambda$-resonant.

\pf Using twice \eqdcomp\ we get
$$
\eqalign{\{\Phi^{-1}\circ F\circ\Phi\}_\nu&=\!\!\!\!\!\!\!
\sum_{{1\le s\le\nu}\atop{\nu_1+\cdots+\nu_s=\nu}}\!\!\!
\tilde K_s(\{F\circ\Phi\}_{\nu_1},\ldots,\{F\circ\Phi\}_{\nu_s})\cr
&=\!\!\!\!\!\!\!\sum_{{1\le s\le\nu}\atop{\nu_1+\cdots+\nu_s=\nu}}
\!\sum_{{1\le r_1\le\nu_1}\atop{d_{11}+\cdots+d_{1r_1}=\nu_1}}\!\!\!\!\!\cdots\!\!\!\!
\sum_{{1\le r_s\le\nu_s}\atop{d_{s1}+\cdots+d_{sr_s}=\nu_s}}\!\!\!\!\!\!\!\!\!\!\!\!\!
\tilde K_s\bigl(\tilde F_{r_1}(H_{d_{11}},\ldots,\!H_{d_{1r_1}}),\ldots,
\tilde F_{r_s}(H_{d_{s1}},\!\ldots,\!H_{d_{sr_s}})\bigr)\cr
&=T_\nu+S_1(\nu)+\sum_{s\ge 2} S_s(\nu)\;,
\cr}
$$
where $\Phi^{-1}=I+\sum\limits_{d\ge 2}K_d$ is the homogeneous expansion
of $\Phi^{-1}$, and:
$$
T_\nu=\sum_{{1\le s\le\nu}\atop{\nu_1+\cdots+\nu_s=\nu}}\!\!\!\tilde K_s(
\Lambda H_{\nu_1},\ldots,\Lambda H_{\nu_s})
\leqno(1)
$$
is obtained considering only the terms with $r_1=\ldots=r_s=1$; 
$$
S_1(\nu)=\sum_{{\mu\le r\le\nu}\atop{d_1+\cdots+d_r=\nu}}\!\!\!
\tilde F_r(H_{d_1},\ldots, H_{d_r})
\leqno (2)
$$
contains the terms with $s=1$ and $r_1>1$; and
$$
S_s(\nu)=\sum_{\nu_1+\cdots+\nu_s=\nu}\!\!
\sum_{{{{1\le r_1\le\nu_1}\atop{\vdots}}\atop{1\le r_s\le\nu_s}}\atop
{\max\{r_1,\ldots,r_s\}\ge\mu}}
\sum_{{{d_{11}+\cdots+d_{1r_1}=\nu_1}\atop{\vdots}}\atop{d_{s1}+\cdots
+d_{sr_s}=\nu_s}}\!\!\!\tilde K_s\bigl(\tilde F_{r_1}(H_{d_{11}},\ldots,
H_{d_{1r_1}}),\ldots,\tilde F_{r_s}(H_{d_{s1}},\ldots,H_{d_{sr_s}})\bigr)
\leqno(3)
$$
contains the terms with fixed $s\ge 2$ and at least one $r_j$ greater than 1
(and thus greater than or equal to $\mu$, because $F_2=\ldots=F_{\mu-1}=O$
by assumption). 

Let us first study $T_\nu$. The summand corresponding to $s=1$ is
$\Lambda H_\nu$; the summand corresponding to $s=\nu$ is $K_\nu\circ\Lambda$; therefore
$$
T_\nu=\Lambda H_\nu+K_\nu\circ\Lambda+\sum_{{2\le s\le\nu-1}\atop{\nu_1+\cdots+\nu_s=\nu}}\!\!\!\tilde K_s(
\Lambda H_{\nu_1},\ldots,\Lambda H_{\nu_s})=-L_\Lambda(H_\nu)+Q_\nu\;,
$$ 
where, using \eqddue\ to express $K_\nu$,
$$
Q_\nu=\sum_{{2\le s\le\nu-1}\atop{\nu_1+\cdots+\nu_s=\nu}}\!\!\!
\left[\tilde K_s(
\Lambda H_{\nu_1},\ldots,\Lambda H_{\nu_s})-
\tilde K_s(H_{\nu_1}\circ\Lambda,\ldots H_{\nu_s}\circ\Lambda)
\right]
$$
depends only on $\Lambda$ and $H_\gamma$ with $\gamma<\nu$ because $2\le s
\le\nu-1$ in the sum. In particular, if $H_1,\ldots,H_{\nu-1}\in\ca H_\Lambda$ then $Q_\nu=O$, and (i) is proved.

Now let us study $S_1(\nu)$. First of all, we clearly have $S_1(\nu)=O$
for $2\le\nu<\mu$, and $S_1(\mu)=F_\mu$. When $\nu>\mu$ we can write
$$
\eqalign{S_1(\nu)&=F_\nu+\sum_{{\mu\le r\le\nu-1}\atop{d_1+\cdots+d_r=\nu}}
\!\!\!\tilde F_r(H_{d_1},\ldots, H_{d_r})\cr
&=F_\nu+\mu\tilde F_\mu(H_{\nu-\mu+1},I,\ldots,I)+\!\!\!\!\!\!\!
\sum_{{d_1+\cdots+d_\mu=
\nu}\atop{1<\max\{d_j\}<\nu-\mu+1}}\!\!\!\!\!\!\tilde F_\mu(H_{d_1},\ldots, H_{d_\mu})
+\!\!\!\sum_{{\mu+1\le r\le\nu-1}\atop{d_1+\cdots+d_r=\nu}}
\!\!\!\tilde F_r(H_{d_1},\ldots, H_{d_r})\;.
\cr}
$$
in particular, $S_1(\mu+1)=F_{\mu+1}+\mu\tilde F_\mu(H_2,I,\ldots,I)$. Notice that the two remaining sums depend only on 
$F_\rho$ with $\rho<\nu$ and on $H_\gamma$ with $\gamma<\nu-\mu+1$
(in the first sum is clear; for the second one, if $d_j\ge\nu-\mu+1$
for some $j$ we then would have $d_1+\cdots+d_r\ge\nu-\mu+1+r-1\ge\nu+1$,  
impossible). Summing up we have
$$
S_1(\nu)=\cases{O&for $2\le\nu<\mu$,\cr
F_\mu&for $\nu=\mu$,\cr
F_{\mu+1}+\mu\tilde F_\mu(H_2,I,\ldots,I)&for $\nu=\mu+1$,\cr
F_\nu+\mu\tilde F_\mu(H_{\nu-\mu+1},I,\ldots,I)+R^1_\nu&for $\nu>\mu+1$,
\cr}
$$
where
$$
R^1_\nu=\!\!\!\sum_{{d_1+\cdots+d_\mu=
\nu}\atop{1<\max\{d_j\}<\nu-\mu+1}}\!\!\!\tilde F_\mu(H_{d_1},\ldots, H_{d_\mu})
+\sum_{{\mu+1\le r\le\nu-1}\atop{d_1+\cdots+d_r=\nu}}
\!\!\!\tilde F_r(H_{d_1},\ldots, H_{d_r})
$$
depends only on 
$F_\rho$ with $\rho<\nu$ and on $H_\gamma$ with $\gamma<\nu-\mu+1$. 

Let us now discuss $S_s(\nu)$ for $s\ge 2$. First of all, the condition
$\max\{r_1,\ldots,r_s\}\ge\mu$ implies
$$
\mu+s-1\le r_1+\cdots+r_s\le\nu_1+\cdots+\nu_s=\nu\;,
$$
that is $s\le\nu-\mu+1$. In particular, $S_s(\nu)=O$ if $\nu\le\mu$ or if
$s>\nu-\mu+1$. 
Moreover, if we had $d_{ij}\ge\nu-\mu+1$ for some $1\le i\le
s$ and $1\le j\le r_s$ we would get
$$
\nu=d_{11}+\cdots+d_{sr_s}\ge\nu-\mu+1+r_1+\cdots+r_s-1
\ge\nu-\mu+1+\mu+s-1-1=\nu+s-1>\nu\;,
$$
impossible. This means that $S_s(\nu)$ depends only on $F_\rho$ with 
$\rho<\nu$ for all $s$, on $H_\gamma$ with $\gamma<\nu-\mu+1$ when $s<\nu-\mu+1$,
and that $S_{\nu-\mu+1}(\nu)$ depends on $H_{\nu-\mu+1}$ just because
it contains $\tilde K_{\nu-\mu+1}$. Furthermore, the
conditions $\max\{r_1,\ldots,r_{\nu-\mu+1}\}\ge\mu$ and $\nu_1+\ldots+\nu_{\nu-\mu+1}=\nu$ imply
that
$$
S_{\nu-\mu+1}(\nu)=(\nu-\mu+1)\tilde K_{\nu-\mu+1}(F_\mu,\Lambda,\ldots,
\Lambda)=-(\nu-\mu+1)\tilde H_{\nu-\mu+1}(F_\mu,\Lambda,\ldots,\Lambda)+R^2_\nu\;,
$$
where (using Lemmas~\udueb\ and~\eqddue)
$$
R^2_\nu=\sum_{{2\le r\le\nu-\mu}\atop{d_1+\cdots+d_r=\nu-\mu+1}}
\sum_{j=1}^r d_j\tilde K_r\bigl(H_{d_1}\circ\Lambda,\ldots,\tilde H_{d_j}(F_\mu,\Lambda,\ldots,
\Lambda),\ldots,H_{d_r}\circ\Lambda\bigr)
$$ 
depends only on $\Lambda$, $F_\mu$ and $H_\gamma$ with
$\gamma<\nu-\mu+1$.

Putting everything together, we have
$$
\eqalign{
\{\Phi^{-1}\circ F\circ\Phi\}_\nu&=T_\nu+S_1(\nu)+\sum_{s=2}^{\nu-\mu+1}
S_s(\nu)\cr
&=F_\nu-L_\Lambda(H_\nu)+Q_\nu+\cases{O&if $2\le\nu\le\mu$,\cr
-L_{F_\mu,\Lambda}(H_2)&if $\nu=\mu+1$,\cr
-L_{F_\mu,\Lambda}(H_{\nu-\mu+1})+R_\nu&if $\nu>\mu+1$,\cr}
\cr}
$$
where
$$
R_\nu=R^1_\nu+R^2_\nu+\sum_{s=2}^{\nu-\mu}S_s(\nu)
$$
depends only on $F_\rho$ with $\rho<\mu$ and on $H_\gamma$ with $\gamma<
\nu-\mu+1$. In particular, if $F=\Lambda$ then we have $S_s(\nu)=O$ for all $s\ge 1$ and hence
$R_\nu=O$ for all $\nu\ge 2$. 

In this way we have proved \eqdqua\ and parts (i), (ii) and (iii). 
Concerning (iv), it suffices to notice that if $F_2,\ldots,F_{\nu-1}$ and $H_2,\ldots, H_{\nu-\mu+1}$
are $\Lambda$-resonant, then also $R^1_\nu$, $S_2(\nu),\ldots,S_{\nu-\mu}(\nu)$ and
$R^2_\nu$ (by Lemmas~\utre\ and~\ddue) are $\Lambda$-resonant.\qedn

\newrem In [WZP2] the remainder term $R_\nu$ is expressed by using combinations of higher order derivatives
instead of combinations of multilinear maps.

We can now introduce the second order normal forms, using the Fischer Hermitian product
to provide suitable complementary spaces.

\newdef The {\sl Fischer Hermitian product} on $\ca H$ is defined by
$$
\langle z_1^{p_1}\cdots z_n^{p_n}e_h,z_1^{q_1}\cdots z_n^{q_n}e_k\rangle=
\cases{
0&if $h\ne k$ or $p_j\ne q_j$ for some $j$;\cr
\noalign{\smallskip}
\displaystyle {p_1!\cdots p_n!\over(p_1+\cdots+p_n)!}&if $h=k$ and $p_j=q_j$ for all $j$.\cr}
$$

\newdef Given $\Lambda\in M_{n,n}(\C)$, we shall say that a $G\in\widehat\ca O^n_\Lambda$ is 
{\sl in second order normal form} if $G=\Lambda$ or the homogeneous expansion $G=\Lambda+\sum\limits_{d\ge\mu}G_d$ of $G$ satisfies the following 
conditions:
\smallskip
\item{(a)} $G_\mu\ne O$; 
\item{(b)} $G_d\in\ca H^d
\cap(\Im L_{G_\mu,\Lambda})^\perp$ for all
$d>\mu$ (where we are using Fischer Hermitian product). 
\smallskip\noindent
Given $F\in\widehat\ca O^n_\Lambda$, we shall say that $G\in\widehat\ca O^n
_\Lambda$ is {\sl a second order normal form} of $F$ if 
$G$ is in second order normal form and $G=\Phi^{-1}\circ F\circ\Phi$ for some $\Phi\in\widehat\ca O^n_I$. 


We can now prove the existence of second order normal forms:

\newthm Theorem \dsei: Let $\Lambda\in M_{n,n}(\C)$ be given. Then each
$F\in\widehat\ca O^n_\Lambda$ admits a second order
normal form. More precisely, if $F=\Lambda+\sum\limits_{d\ge\mu}F_d$ is in Poincar\'e-Dulac normal form (and $F\not\equiv\Lambda$) then there exists a unique 
$\Lambda$-resonant $\Phi=I+\sum\limits_{d\ge 2}H_d\in\widehat\ca O^n_I$ such that $H_d\in(\Ker L_{F_\mu,\Lambda})^\perp$ for all $d\ge 2$ and $G=\Phi^{-1}\circ F\circ\Phi$ is in second order
normal form. Furthermore, if $\Lambda$ is diagonal we also have $G_d\in\ca H^d_\Lambda$ for all $d\ge\mu$. 

\pf By the classical theory we can assume that $F$ is in Poincar\'e-Dulac
normal form. If $F\equiv\Lambda$ we are done; assume then that $F\not\equiv\Lambda$. 

First of all, by Theorem~\dquattro\ if $\Phi$ is $\Lambda$-resonant 
we have $\{\Phi^{-1}\circ F\circ\Phi\}_d=F_d$ for all $d\le\mu$.  Now consider the splittings
$$
\ca H^d=\Im L_{F_\mu,\Lambda}|_{\ca H^{d-\mu+1}_\Lambda}\bigoperp
(\Im L_{F_\mu,\Lambda}|_{\ca H^{d-\mu+1}_\Lambda})^\perp
$$
and 
$$
\ca H^{d-\mu+1}_\Lambda=\Ker L_{F_\mu,\Lambda}|_{\ca H^{d-\mu+1}_\Lambda}\bigoperp
(\Ker L_{F_\mu,\Lambda}|_{\ca H^{d-\mu+1}_\Lambda})^\perp\;.
$$
If $d=\mu+1$ we can find a unique $G_{\mu+1}\in(\Im L_{F_\mu,\Lambda})^\perp
\cap\ca H^{\mu+1}$ and a unique $H_2\in (\Ker L_{F_\mu,\Lambda})^\perp
\cap\ca H^{2}_\Lambda$ such that $F_{\mu+1}=G_{\mu+1}+L_{F_\mu,\Lambda}(H_2)$. Then 
Theorem~\dquattro\ yields
$$
\{\Phi^{-1}\circ F\circ\Phi\}_{\mu+1}=F_{\mu+1}-L_{F_\mu,\Lambda}(\{\Phi\}_2)
=G_{\mu+1}+L_{F_\mu,\Lambda}(H_2)-L_{F_\mu,\Lambda}(\{\Phi\}_2)\;;
$$
so to get $\{\Phi^{-1}\circ F\circ\Phi\}_{\mu+1}\in(\Im L_{F_\mu,\Lambda})^\perp
\cap\ca H^{\mu+1}$ with $\{\Phi\}_2\in (\Ker L_{F_\mu,\Lambda})^\perp
\cap\ca H^{2}_\Lambda$ we must necessarily take $\{\Phi\}_2=H_2$.

Assume, by induction, that we have uniquely determined $H_2,\ldots,H_{d-\mu}\in
(\Im L_{F_\mu,\Lambda})^\perp\cap\ca H_\Lambda$, and thus
$R_d\in\ca H^d$ in \eqdqua. Hence there is a unique $G_d\in(\Im L_{F_\mu,\Lambda})^\perp
\cap\ca H^d$ and a unique $H_{d-\mu+1}\in (\Ker L_{F_\mu,\Lambda})^\perp
\cap\ca H^{d-\mu+1}_\Lambda$ such that $F_d+R_d=G_d+L_{F_\mu,\Lambda}(H_{d-\mu+1})$. Thus
to get $\{\Phi^{-1}\circ F\circ\Phi\}_d\in(\Im L_{F_\mu,\Lambda})^\perp
\cap\ca H^d$ with $\{\Phi\}_{d-\mu+1}\in(\Ker L_{F_\mu,\Lambda})^\perp
\cap\ca H^{d-\mu+1}_\Lambda$ the 
only possible choice is $\{\Phi\}_{d-\mu+1}=H_{d-\mu+1}$, and thus $\{\Phi^{-1}\circ F\circ\Phi\}_d=G_d$.

Finally, if $\Lambda$ is diagonal then $F_d\in\ca H^d_\Lambda$ for all $d\ge\mu$. 
Furthermore, Lemma~\uqua\ yields $\Im L_{F_\mu,\Lambda}|_{\ca H^{d-\mu+1}_\Lambda}
\subseteq\ca H^d_\Lambda$ for all $d\ge\mu$; recalling Theorem~\dquattro.(vi)
we then see can we can always find $G_d\in\ca H_\Lambda^d$, and we are done.\qedn

The definition and construction of $k$-th order normal forms is similar; the idea is to extract from
the remainder term $R_\nu$ the pieces depending on $H_\gamma$ with $\gamma$ varying in a suitable range, and use them to build operators generalizing $L_\Lambda$ and $L_{P,\Lambda}$. 
We refer to [WZP2] for details; for our needs it suffices to recall that given $F=\Lambda+\sum_{d\ge 2} F_d\in\widehat{\ca O}^n_\Lambda$ [WPZ2] introduces a sequence of (not necessarily linear)
operators $\ca L^{(d)}[\Lambda,F_2,\ldots,F_d]\colon\Ker\ca L^{(d-1)}\times\ca H^{d+1}\to\ca H^{d+1}$ for $d\ge 1$, with $\ca L^{(1)}[\Lambda](H_2)=L_\Lambda(H_2)$ and
$\ca L^{(2)}[\Lambda, F_2](H_2,H_3)=L_\Lambda(H_3)+L_{F_2,\Lambda}(H_2)$, and gives
the following definition:

\newdef We shall say that $G=\Lambda+\sum\limits_{d\ge 2} G_d\in\widehat{\ca O}^n_\Lambda$
is in {\sl infinite order normal form} if $G_d\in W_d^\perp$ for all $d\ge 2$, where $W_d$
is a vector subspace of maximal dimension contained in the image of $\ca L^{(d-1)}[\Lambda,G_2,\ldots,G_{d-1}]$. We shall also say that
$G$ is an {\sl infinite order normal form of $F\in\widehat{\ca O}^n_\Lambda$} if it is in infinite order normal form and
it is formally conjugated to~$F$.

We end this section quoting a result from [WZP2] giving a condition ensuring that
a second order normal form is actually an infinite order normal form:

\newthm Proposition \dinf: ([WZP2, Theorem~4.9]) Let $\Lambda\in M_{n,n}(\C)$ be diagonal,
and $F=\Lambda+\sum\limits_{d\ge 2}F_d\in\widehat{\ca O}^n_\Lambda$ with $F_2\ne O$
and $\Lambda$-resonant. Assume that $\Ker L_{F_2,\Lambda}|_{\ca H^d_\Lambda}=\{O\}$
for all $d\ge 2$. Then the second order normal form of~$F$ is the unique infinite order normal
form of~$F$.

\smallsect 2. Superattracting germs

In this section we shall completely describe the second order normal forms obtained when $n=\mu=2$ and $\Lambda=O$, that is for 2-dimensional superattracting germs with non-vanishing quadratic term. Except in four degenerate
instances, the second order normal forms will be infinite order normal forms, and will be expressed just in terms of two power series of {\it one} variable, thus giving a drastic simplification of the germs.

In [A1] we 
showed that, up to a linear change of variable, we can assume that the quadratic term $F_2$ 
is of one (and only one) of the following forms:
\smallskip
\itemitem{$(\infty)$} $F_2(z,w)=(z^2,zw)$;
\itemitem{$(1_{00})$} $F_2(z,w)=(0,-z^2)$;
\itemitem{$(1_{10})$} $F_2(z,w)=\bigl(-z^2,-(z^2+zw)\bigr)$;
\itemitem{$(1_{11})$} $F_2(z,w)=\bigl(-zw,-(z^2+w^2)\bigr)$;
\itemitem{$(2_{001})$} $F_2(z,w)=(0,zw)$;
\itemitem{$(2_{011})$} $F_2(z,w)=(zw,zw+w^2)$;
\itemitem{$(2_{10\rho})$} $F_2(z,w)=\bigl(-\rho z^2,(1-\rho) zw)$, with $\rho\ne0$;
\itemitem{$(2_{11\rho})$} $F_2(z,w)=\bigl(\rho z^2+zw,(1+\rho) zw+w^2\bigr)$, with $\rho\ne0$;
\itemitem{$(3_{100})$} $F_2(z,w)=(z^2-zw,0)$;
\itemitem{$(3_{\rho10})$} $F_2(z,w)=\bigl(\rho(-z^2+zw),(1-\rho)(zw- w^2)\bigr)$, with
$\rho\ne0$,~$1$;
\itemitem{$(3_{\rho\tau1})$} $F_2(z,w)=\bigl(-\rho z^2+(1-\tau)zw,(1-\rho)
zw-\tau w^2\bigr)$, with $\rho$,~$\tau\ne0$ and $\rho+\tau\ne1$
\smallskip
\noindent (where the symbols refer to the number of characteristic directions and to their indeces;
see also [AT2]). 


We shall use the standard basis $\{u_{d,j},v_{d,j}\}_{j=0,\ldots,d}$ of $\ca H^d$, where
$$
u_{d,j}=(z^jw^{d-j},0)\qquad\hbox{and}\qquad v_{d,j}=(0,z^jw^{d-j})\;,
$$
and we shall endow $\ca H^d$ with  Fischer Hermitian product, so that $\{u_{d,j},v_{d,j}\}_{j=0,\ldots,d}$ is an orthogonal basis and
$$
\|u_{d,j}\|^2=\|v_{d,j}\|^2={d\choose j}^{-1}\;.
$$

When $\Lambda=O$, we have $\ca H_\Lambda=\ca H$, and the operator $L=L_{F_2,\Lambda}\colon\ca H^d\to\ca H^{d+1}$ is given by
$$
L(H)=-\Jac(F_2)\cdot H\;.
$$
To apply Proposition~\dinf, we need to know when $\Ker L|_{\ca H^d}=\{O\}$. Since
$$
\dim\Ker L|_{\ca H^d}+\dim\Im L|_{\ca H^d}=\dim\ca H^d=\dim\ca H^{d+1}-2
=\dim\Im L|_{\ca H^d}+\dim(\Im L|_{\ca H^d})^\perp-2\;,
$$
we find that
$$
\Ker L|_{\ca H^d}=\{O\}\quad\hbox{if and only if}\quad \dim(\Im L|_{\ca H^d})^\perp=2\;.
\neweq\eqdinfo
$$

We shall now study separately each case.

\medbreak
\noindent$\bullet$ {\it Case $(\infty)$.}

\noindent In this case we have
$$
L(u_{d,j})=-2u_{d+1,j+1}-v_{d+1,j}\qquad\hbox{and}\qquad L(v_{d,j})=-v_{d+1,j+1}
$$
for all $d\ge 2$ and $j=0,\ldots,d$. Therefore
$$
\Im L|_{\ca H^d}=\Span{u_{d+1,2},\ldots,u_{d+1,d+1}, 2u_{d+1,1}+v_{d+1,0},v_{d+1,1},\ldots,v_{d+1,d+1}}\;,
$$
and thus
$$
(\Im L|_{\ca H^d})^\perp=\Span{u_{d+1,0},(d+1)u_{d+1,1}-2v_{d+1,0}}\;.
$$
In particular, thanks to \eqdinfo\ and Proposition~\dinf, a second order normal form is automatically an infinite
order normal form.

It then follows that every formal power series of the form
$$
F(z,w)=(z^2+O_3,zw+O_3)
$$
(where $O_3$ denotes a remainder term of order at least~3)
has a unique infinite order normal form 
$$
G(z,w)=\bigl(z^2+\phe(w)+z\psi'(w),zw-2\psi(w)\bigr)
$$
where $\phe$, $\psi\in\C[\![\zeta]\!]$ are power series of order at least~3. Notice that
(here and in later formulas) the appearance of the derivative (which simplifies the expression
of the normal form) is due to the fact we are using the Fischer Hermitian product; using
another Hermitian product might lead to more complicated normal
forms. 

\medbreak
\noindent$\bullet$ {\it Case $(1_{00})$.}

\noindent In this case we have
$$
L(u_{d,j})=2v_{d+1,j+1}\qquad\hbox{and}\qquad L(v_{d,j})=0
$$
for all $d\ge 2$ and $j=0,\ldots,d$. Therefore
$$
\Im L|_{\ca H^d}=\Span{v_{d+1,1},\ldots,v_{d+1,d+1}}\;,
$$
and thus
$$
(\Im L|_{\ca H^d})^\perp=\Span{u_{d+1,0},\ldots,u_{d+1,d+1},v_{d+1,0}}\;.
$$
This a degenerate case, where we cannot use Proposition~\dinf. Anyway, Theorem~\dsei\
still apply, and 
it follows that every formal power series of the form
$$
F(z,w)=(O_3,-z^2+O_3)
$$
has a second order normal form
$$
G(z,w)=\left(\Phi(z,w),-z^2+\psi(w)\right)
$$
where $\psi\in\C[\![\zeta]\!]$ and $\Phi\in\C[\![z,w]\!]$ are power series of
order at least~3.

\medbreak
\noindent$\bullet$ {\it Case $(1_{10})$.}

\noindent In this case we have
$$
L(u_{d,j})=2u_{d+1,j+1}+2v_{d+1,j+1}+v_{d+1,j}\qquad\hbox{and}\qquad L(v_{d,j})=v_{d+1,j+1}
$$
for all $d\ge 2$ and $j=0,\ldots,d$. Therefore
$$
\Im L|_{\ca H^d}=\Span{2u_{d+1,1}+v_{d+1,0},u_{d+1,2},\ldots,u_{d+1,d+1},v_{d+1,1},\ldots,v_{d+1,d+1}}\;,
$$
and thus
$$
(\Im L|_{\ca H^d})^\perp=\Span{u_{d+1,0},(d+1)u_{d+1,1}-2v_{d+1,0}}\;.
$$
It then follows that every formal power series of the form
$$
F(z,w)=(-z^2+O_3,-z^2-zw+O_3)
$$
has a unique infinite order normal form
$$
G(z,w)=\bigl(-z^2+\phe(w)+z\psi'(w),-z^2-zw-2\psi(w)\bigr)
$$
where $\phe$,~$\psi\in\C[\![\zeta]\!]$ are power series of order at least~3.

\medbreak
\noindent$\bullet$ {\it Case $(1_{11})$.}

\noindent In this case we have
$$
L(u_{d,j})=u_{d+1,j}+2v_{d+1,j+1}\qquad\hbox{and}\qquad L(v_{d,j})=u_{d+1,j+1}+2v_{d+1,j}
$$
for all $d\ge 2$ and $j=0,\ldots,d$. It follows that 
$$
\eqalign{
\Im L|_{\ca H^d}=\hbox{\rm Span}(&u_{d+1,0}-u_{d+1,2},\ldots,u_{d+1,d-1}-u_{d+1,d+1},\cr
&\quad v_{d+1,2}-v_{d+1,0},
\ldots,v_{d+1,d+1}-v_{d+1,d-1},u_{d+1,0}+2v_{d+1,1},u_{d+1,1}+2v_{d+1,0})\;,\cr}
$$
and a few computations yield
$$
\eqalign{
(\Im L|_{\ca H^d})^\perp&=\Span{\sum_{j=0}^{d+1}{d+1\choose j}(v_{d+1,j}-2u_{d+1,j}),
\sum_{j=0}^{d+1}(-1)^j{d+1\choose j}(v_{d+1,j}+2u_{d+1,j})}\cr
&=\Span{\bigl(-2(z+w)^{d+1},(z+w)^{d+1}\bigr),\bigl(2(w-z)^{d+1},(w-z)^{d+1}\bigr)}\cr}\;.
$$
It then follows that every formal germ of the form
$$
F(z,w)=(-zw+O_3,-z^2-w^2+O_3)
$$
has a unique infinite order normal form
$$
G(z,w)=\bigl(-zw-2\phe(z+w)+2\psi(w-z),
-z^2-w^2+\phe(z+w)+
\psi(w-z)\bigr)
$$
where $\phe$,~$\psi\in\C[\![\zeta]\!]$ are arbitrary power series of order at least~3. Again, the fact that the normal form is expressed in terms of power series evaluated in $z+w$ and $z-w$
is due to the fact we are using Fischer Hermitian product.

\medbreak
\noindent$\bullet$ {\it Case $(2_{001})$.}

\noindent In this case we have
$$
L(u_{d,j})=-v_{d+1,j}\qquad\hbox{and}\qquad L(v_{d,j})=-v_{d+1,j+1}
$$
for all $d\ge 2$ and $j=0,\ldots,d$. It follows that 
$$
\Im L|_{\ca H^d}=\Span{v_{d+1,0},\ldots,v_{d+1,d+1}}
$$
and hence
$$
(\Im L|_{\ca H^d})^\perp=\Span{u_{d+1,0},\ldots,u_{d+1,d+1}}\;.
$$
We are in a degenerate case; hence
every formal germ of the form
$$
F(z,w)=(O_3,zw+O_3)
$$
has a second order normal form
$$
G(z,w)=\bigl(\Phi(z,w),zw\bigr)
$$
where $\Phi\in\C[\![z,w]\!]$ is a power series of order at least three.

\medbreak
\noindent$\bullet$ {\it Case $(2_{011})$.}

\noindent In this case we have
$$
L(u_{d,j})=-u_{d+1,j}-v_{d+1,j}\qquad\hbox{and}\qquad L(v_{d,j})=-u_{d+1,j+1}-2v_{d+1,j}
-v_{d+1,j+1}
$$
for all $d\ge 2$ and $j=0,\ldots,d$. It follows that 
$$
\eqalign{
\Im L|_{\ca H^d}&=\hbox{\rm Span}(u_{d+1,0},\ldots,u_{d+1,d-1},v_{d+1,0},\ldots,v_{d+1,d-1},\cr
&\qquad\qquad u_{d+1,d}+v_{d+1,d},u_{d+1,d+1}+v_{d+1,d+1}+2v_{d+1,d})\;,\cr}
$$
and hence
$$
(\Im L|_{\ca H^d})^\perp=\Span{(d+1) u_{d+1,d}- (d+1) v_{d+1,d}+2v_{d+1,d+1},u_{d+1,d+1}-v_{d+1,d+1}}\;.
$$
It then follows that every formal germ of the form
$$
F(z,w)=(zw+O_3,zw+w^2+O_3)
$$
has a unique infinite order normal form
$$
G(z,w)=\bigl( zw + w\phe'(z) + \psi(z),
zw + w^2 + 2\phe(z) - w\phe'(z) - \psi(z)\bigr)\;,
$$
where $\phe$, $\psi\in\C[\![\zeta]\!]$ are power series of order at least~3.

\medbreak
\noindent$\bullet$ {\it Case $(2_{10\rho})$.}

\noindent In this case we have
$$
L(u_{d,j})=2\rho u_{d+1,j+1}+(\rho-1)v_{d+1,j}\qquad\hbox{and}\qquad L(v_{d,j})=(\rho-1)v_{d+1,j+1}$$
for all $d\ge 2$ and $j=0,\ldots,d$. We clearly have two subcases to consider.

If $\rho=1$ then
$$
\Im L|_{\ca H^d}=\Span{u_{d+1,1},\ldots,u_{d+1,d+1}}\;,
$$
and hence
$$
(\Im L|_{\ca H^d})^\perp=\Span{u_{d+1,0},v_{d+1,0},\ldots,v_{d+1,d+1}}\;.
$$
We are in the third degenerate case; hence every formal germ of the form
$$
F(z,w)=(-z^2+O_3,O_3)
$$
has a second order normal form
$$
G(z,w)=\bigl(-z^2+\psi(w),\Phi(z,w)\bigr)\;,
$$
where $\psi\in\C[\![\zeta]\!]$ and $\Phi\in\C[\![z,w]\!]$ are power series of order at least~3.

If instead $\rho\ne 1$ (recalling that $\rho\ne0$ too) then
$$
\Im L|_{\ca H^d}=\Span{2\rho u_{d+1,1} + (\rho - 1) v_{d+1,0}, u_{d+1,2},\ldots,u_{d+1,d+1},v_{d+1,1},\ldots,v_{d+1,d+1}}\;,
$$
and hence
$$
(\Im L|_{\ca H^d})^\perp=\Span{u_{d+1,0},(\rho-1)(d+1) u_{d+1,1} - 2\rho v_{d+1,0}}\;.
$$
It then follows that every formal germ of the form
$$
F(z,w)=\bigl(-\rho z^2+O_3,(1-\rho)zw+O_3\bigr)
$$
with $\rho\ne 0$,~$1$ has a unique infinite order normal form
$$
G(z,w)=\bigl(- \rho z^2+ (\rho -1) z\phe'(w)+ \psi(w),(1-\rho)zw - 2\rho\phe(z)\bigr)\;,
$$
where $\phe$,~$\psi\in\C[\![\zeta]\!]$ are power series of order at least~3.

\medbreak
\noindent$\bullet$ {\it Case $(2_{11\rho})$.}

\noindent In this case we have
$$
\cases{
L(u_{d,j})=-2\rho u_{d+1,j+1}-u_{d+1,j}-(1+\rho)v_{d+1,j}\cr
L(v_{d,j})=
-u_{d+1,j+1}-2v_{d+1,j}-(1+\rho)v_{d+1,j+1}\cr}
\neweq\eqtuno
$$
for all $d\ge 2$ and $j=0,\ldots,d$. We clearly have two subcases to consider.

If $\rho=-1$ then
$$
\Im L|_{\ca H^d}=\Span{u_{d+1,0}-2u_{d+1,1},\ldots,u_{d+1,d}-2u_{d+1,d+1},
u_{d+1,1}+2v_{d+1,0},\ldots,u_{d+1,d}+2v_{d+1,d}}\;,
$$
and hence
$$
\eqalign{
(\Im L|_{\ca H^d})^\perp&=\Span{\sum_{j=0}^{d+1}{d+1\choose j}{1\over 2^j}(u_{d+1,j}-{1\over 4}
v_{d+1,j}),v_{d+1,d+1}}\cr
&=\Span{\left(\left({\textstyle{{z\over 2}}}+w\right)^{d+1},-{1\over 4}\left({\textstyle{{z\over 2}}}+w\right)^{d+1}
\right),(0,z^{d+1})}\;.\cr}
$$
It then follows that every formal germ of the form
$$
F(z,w)=(-z^2+zw+O_3,w^2+O_3)
$$
has a unique infinite order normal form
$$
G(z,w)=\left(-z^2+zw+\phe({\textstyle{{z\over 2}}}+w),
w^2- {1\over 4}\phe({\textstyle{{z\over 2}}}+w)+\psi(z)\right)\;,
$$
where $\phe$,~$\psi\in\C[\![\zeta]\!]$ are power series of order at least~3.

If instead $\rho\ne -1$ (recalling that $\rho\ne0$ too) then a basis of $\Im L|_{\ca H^d}$
is given by the vectors listed in \eqtuno, and a computation shows that $(\Im L|_{\ca H^d})^\perp$
is given by homogeneous maps of the form
$$
\sum_{j=0}^{d+1}(a_j u_{d+1,j}+b_j v_{d+1,j})
$$
where the coefficients $a_j$, $b_j$ satisfy the following relations:
$$
\cases{\displaystyle c_jb_j=-{2\over1+\rho}c_{j-1}b_{j-1}-{1\over\rho(1+\rho)}c_{j-2}b_{j-2}&for $j=2,\ldots,d+1$,\cr
\noalign{\smallskip}
\displaystyle c_ja_j={1\over\rho}c_{j-2}b_{j-2}&for $j=2,\ldots,d+1$,\cr
\noalign{\smallskip}
\displaystyle a_0=(3\rho-1)b_0+2{\rho(1+\rho)\over d+1}b_1\;,\cr
\noalign{\smallskip}
\displaystyle a_1=-2(d+1)b_0-(1+\rho)b_1\;,\cr}
$$
where $c_j^{-1}={d+1\choose j}$
and $b_0$, $b_1\in\C$ are arbitrary. Solving these recurrence equations one gets
$$
b_j={1\over2\sqrt{-\rho}}{d+1\choose j}\left[{\rho(1+\rho)\over d+1}(m_\rho^j-n_\rho^j)b_1
+\bigl(\rho(m_\rho^j-n_\rho^j)+\sqrt{-\rho}(m_\rho^j+n_\rho^j)\bigr)b_0\right]\;,
$$
where $\sqrt{-\rho}$ is any square root of $-\rho$, and
$$
m_\rho={\sqrt{-\rho}-\rho\over\rho(1+\rho)}\;,\qquad n_\rho=-{\sqrt{-\rho}+\rho\over\rho(1+\rho)}\;.
$$
It follows that the unique infinite order normal form of a formal germ of the form
$$
F(z,w)=\bigl(\rho z^2+zw+O_3,(1+\rho)zw+w^2+O_3\bigr)
$$
with $\rho\ne 0$,~$-1$ is
$$
\eqalign{
G(z,w)=&\left(\rho z^2+zw+{1\over\rho}\left[{1-\sqrt{-\rho}\over 2m_\rho^2}\phe(m_\rho z+w)
+{1+\sqrt{-\rho}\over 2n_\rho^2}\phe(n_\rho z+w)\right]\right.\cr
&\phantom{\biggl(\rho z^2+zw\,}
+{1+\rho\over2\sqrt{-\rho}}\left(
{1\over m_\rho^2}\psi(m_\rho z+w)-{1\over n_\rho^2}\psi(n_\rho z+w)\right),\cr
&\quad(1+\rho)zw+w^2+ {1-\sqrt{-\rho}\over2}\phe(m_\rho z+w)+{1+\sqrt{-\rho}\over2}\phe(n_\rho z+w)\cr
&\phantom{\biggl((1+\rho)zw+w^2\;}+\left.{\rho(1+\rho)\over2\sqrt{-\rho}}\bigl(\psi(m_\rho z+w)-\psi(n_\rho z+w)\bigr)
\right)\cr}
$$ 
where $\phe$,~$\psi\in\C[\![\zeta]\!]$ are power series of order at least~3.

\medbreak
\noindent$\bullet$ {\it Case $(3_{100})$.}

\noindent In this case we have
$$
L(u_{d,j})=u_{d+1,j}-2u_{d+1,j+1}\qquad\hbox{and}\qquad L(v_{d,j})=u_{d+1,j+1}
$$
for all $d\ge 2$ and $j=0,\ldots,d$. It follows that 
$$
\Im L|_{\ca H^d}=\Span{u_{d+1,0},\ldots,u_{d+1,d+1}}
$$
and hence
$$
(\Im L|_{\ca H^d})^\perp=\Span{v_{d+1,0},\ldots,v_{d+1,d+1}}\;.
$$
We are in the last degenerate case; hence every formal germ of the form
$$
F(z,w)=(z^2-zw+O_3,O_3)
$$
has a second order normal form
$$
G(z,w)=\bigl(z^2-zw,\Phi(z,w)\bigr)\;,
$$
where $\Phi\in\C[\![z,w]\!]$ is a power series of order at least~3.

\medbreak
\noindent$\bullet$ {\it Case $(3_{\rho10})$.}

\noindent In this case we have
$$
\cases{
L(u_{d,j})=\rho(2u_{d+1,j+1}-u_{d+1,j})+(\rho-1)v_{d+1,j}\cr
L(v_{d,j})=-\rho u_{d+1,j+1}+(\rho-1)(v_{d+1,j+1}-2v_{d+1,j})\cr}
\neweq\eqtdue
$$
for all $d\ge 2$ and $j=0,\ldots,d$. Then a basis of $\Im L|_{\ca H^d}$ is given by the 
homogeneous maps listed in \eqtdue, and a computation shows that $(\Im L|_{\ca H^d})^\perp$
is given by homogeneous maps of the form
$$
\sum_{j=0}^{d+1}(a_j u_{d+1,j}+b_j v_{d+1,j})
$$
where the coefficients $a_j$, $b_j$ satisfy the following relations:
$$
\cases{\displaystyle c_{j+1}a_{j+1}={\rho-1\over\rho}(c_{j+1}b_{j+1}-2c_jb_j)&for $j=0,\ldots,d$,\cr
\noalign{\smallskip}
\displaystyle c_{j+1}b_{j+1}=2c_jb_j-c_{j-1}b_{j-1}&for $j=1,\ldots,d$,\cr
\noalign{\smallskip}
\displaystyle c_0a_0=2c_1a_1+{\rho-1\over \rho}c_0b_0\;,\cr}
$$
where $c_j^{-1}={d+1\choose j}$ and $b_0$, $b_1\in\C$ are arbitrary. Solving these recurrence
equations we find
$$
\cases{
b_j={d+1\choose j}\left[{j\over d+1}b_1-(j-1)b_0\right]&
for $j=0,\ldots,d+1$,\cr
a_j={\rho-1\over \rho}{d+1\choose j}\left[{2-j\over d+1}b_1+(j-3)b_0
\right]&for $j=0,\ldots,d+1$,\cr}
$$
where $b_0$, $b_1\in\C$ are arbitrary. So every formal germ of the form
$$
F(z,w)=\bigl(\rho (-z^2+zw)+O_3,(1-\rho)(zw-w^2)+O_3\bigr)
$$
with $\rho\ne 0$,~$1$ has a unique infinite order normal form
$$
\eqalign{
G(z,w)=&\left(\rho(-z^2+zw)+z{\de\over\de z}\bigl[\phe(z+w)+\psi(z+w)\bigr]-\phe(z+w),\right.\cr
&\left.\quad(1-\rho)(zw-w^2)+{\rho-1\over\rho}\left(z{\de\over\de z}\bigl[\phe(z+w)-\psi(z+w)\bigr]
-3\phe(z+w)+2\psi(z+w)\right)
\right)\cr}
$$ 
where $\phe$,~$\psi\in\C[\![\zeta]\!]$ are power series of order at least~3.

\medbreak
\noindent$\bullet$ {\it Case $(3_{\rho\tau1})$.}

\noindent In this case we have
$$
L(u_{d,j})=(\tau-1)u_{d+1,j}+2\rho u_{d+1,j+1}+(\rho-1)v_{d+1,j}
$$
and
$$  
L(v_{d,j})=(\tau-1)u_{d+1,j+1}+2\tau v_{d+1,j}+(\rho-1)v_{d+1,j+1}
$$
for all $d\ge 2$ and $j=0,\ldots,d$. As before, we have a few subcases to consider.

Assume first $\rho=\tau=1$. Then
$$
\Im L|_{\ca H^d}=\Span{u_{d+1,1},\ldots,u_{d+1,d+1},v_{d+1,0},\ldots,v_{d+1,d}}\;;
$$
hence
$$
(\Im L|_{\ca H^d})^\perp=\Span{u_{d+1,0},v_{d+1,d+1}}\;,
$$
It then follows that every formal germ of the form
$$
F(z,w)=(-z^2+O_3,-w^2+O_3)
$$
has a unique infinite order normal form
$$
G(z,w)=\bigl(-z^2+\phe(w),-w^2+\psi(z)\bigr)\;,
$$
where $\phe$,~$\psi\in\C[\![\zeta]\!]$ are power series of order at least~3.

Assume now $\rho\ne 1$. Then a computation shows that $(\Im L|_{\ca H^d})^\perp$
is given by homogeneous maps of the form
$$
\sum_{j=0}^{d+1}(a_j u_{d+1,j}+b_j v_{d+1,j})
$$
where the coefficients $a_j$, $b_j$ satisfy the following relations:
$$
\cases{\displaystyle c_{j+1}a_{j+1}={\tau\over\rho}c_{j-1}b_{j-1}&for $j=1,\ldots,d$,\cr
\noalign{\smallskip}
\displaystyle c_{j+1}b_{j+1}=-{2\tau\over\rho-1}c_jb_j-{\tau(\tau-1)\over\rho(\rho-1)}c_{j-1}b_{j-1}&for $j=1,\ldots,d$,\cr
\noalign{\smallskip}
\displaystyle (\tau-1)c_1a_1+(\rho-1)c_1b_1+2\tau c_0b_0=0\;,\cr
\displaystyle (\tau-1)c_0a_0+(\rho-1)c_0b_0+2\rho c_1a_1=0\;,\cr}
\neweq\eqttre
$$
where $c_j^{-1}={d+1\choose j}$ and $b_0$, $b_1\in\C$ are arbitrary. 

When $\tau=1$ conditions \eqttre\ reduce to
$$
\cases{\displaystyle c_{j+1}a_{j+1}={1\over\rho}c_{j-1}b_{j-1}&for $j=1,\ldots,d$,\cr
\noalign{\smallskip}
\displaystyle c_{j+1}b_{j+1}=-{2\over\rho-1}c_jb_j&for $j=1,\ldots,d$,\cr
\noalign{\smallskip}
\displaystyle (\rho-1)c_1b_1+2 c_0b_0=0\;,\cr
\displaystyle (\rho-1)c_0b_0+2\rho c_1a_1=0\;,\cr}
$$
whose solution is
$$
\cases{
\displaystyle a_j={d+1\choose j}{1\over\rho}\left({2\over 1-\rho}\right)^{j-2}b_0&for $j=1,\ldots,d+1$,\cr
\noalign{\smallskip}
\displaystyle b_j={d+1\choose j}\left({2\over1-\rho}\right)^j b_0&for $j=0,\ldots, d+1$,\cr}
$$
where $a_0$, $b_0\in\C$ are arbitrary. Therefore
$$
(\Im L|_{\ca H^d})^\perp=\Span{(w^{d+1},0),\left({(1-\rho)^2\over4\rho}\left({2\over1-\rho}z+w
\right)^{d+1}, \left({2\over1-\rho}z+w\right)^{d+1}\right)}\;,
$$
and thus every formal germ of the form
$$
F(z,w)=\bigl(-\rho z^2+O_3,(1-\rho)zw- w^2+O_3\bigr)
$$
with $\rho\ne 1$ has a unique infinite order normal form
$$
G(z,w)=\left(-\rho z^2+\phe(w)+{(1-\rho)^2\over4\rho}\psi\left({2\over1-\rho}z+w\right),
(1-\rho)zw-w^2+\psi\left(
{2\over1-\rho}z+w\right)\right)\;,
$$
where $\phe$,~$\psi\in\C[\![\zeta]\!]$ are arbitrary power series of order at least~3.

The case $\rho=1$ and $\tau\ne 1$ is treated in the same way; we get that
every formal germ of the form
$$
F(z,w)=\bigl(- z^2+(1-\tau)zw+O_3,-\tau w^2+O_3\bigr)
$$
with $\tau\ne 1$ has a unique infinite order normal form
$$
G(z,w)=\left(- z^2+(1-\tau)zw+\psi\left({1-\tau\over 2}z+w\right),
-\tau w^2+\phe(z)+{(1-\tau)^2\over 4\tau}\psi\left(
{1-\tau\over 2}z+w\right)\right)\;,
$$
where $\phe$,~$\psi\in\C[\![\zeta]\!]$ are power series of order at least~3.

Finally assume $\rho$,~$\tau\ne 1$ (and $\rho+\tau\ne 1$). Solving the recurrence
equations \eqttre\ we find
$$
\eqalign{
b_j={1\over2\sqrt{\rho\tau(\rho+\tau-1)}}{d+1\choose j}&\left[{\rho(\rho-1)\over d+1}(m_{\rho,\tau}^j-n_{\rho,\tau}^j)b_1\right.\cr
&\quad+\bigl(\rho\tau(m_{\rho,\tau}^j-n_{\rho,\tau}^j)+\sqrt{\rho\tau(\rho+\tau-1)}(m_{\rho,\tau}^j+n_{\rho,\tau}^j)\bigr)b_0\biggr]\;,\cr}
$$
for $j=0,\ldots,d+1$,
where $\sqrt{\rho\tau(\rho+\tau-1)}$ is any square root of $\rho\tau(\rho+\tau-1)$, and
$$
m_{\rho,\tau}={\sqrt{\rho\tau(\rho+\tau-1)}-\rho\tau\over\rho(\rho-1)}\;,\qquad n_{\rho,\tau}
=-{\sqrt{\rho\tau(\rho+\tau-1)}+\rho\tau\over\rho(\rho-1)}\;.
$$
Moreover, from \eqttre\ we also get
$$
\eqalign{
a_j={\tau\over2\rho\sqrt{\rho\tau(\rho+\tau-1)}}{d+1\choose j}&\left[{\rho(\rho-1)\over d+1}(m_{\rho,\tau}^{j-2}-n_{\rho,\tau}^{j-2})b_1\right.\cr
&\quad+\bigl(\rho\tau(m_{\rho,\tau}^{j-2}-n_{\rho,\tau}^{j-2})+\sqrt{\rho\tau(\rho+\tau-1)}(m_{\rho,\tau}^{j-2}+n_{\rho,\tau}^{j-2})\bigr)b_0\biggr]\;,\cr}
$$
again for $j=0,\ldots,d+1$.
It follows that the unique infinite order normal form of a formal germ of the form
$$
F(z,w)=\bigl(-\rho z^2+(1-\tau)zw+O_3,(1-\rho)zw-\tau w^2+O_3\bigr)
$$
with $\rho$, $\tau\ne 0$,~1 and $\rho+\tau\ne 1$, is
$$
\eqalign{
G(z,w)=&\left(-\rho z^2+(1-\tau)zw
+{\tau\over\rho}\left[{\sqrt{\rho+\tau-1}+\sqrt{\rho\tau}\over 2m_{\rho,\tau}^2}\phe(m_{\rho,\tau}z+w)
\right.\right.\cr
&\phantom{-\rho z^2+(1-\tau)zw+{\tau\over\rho}\biggl[\quad}+{\sqrt{\rho+\tau-1}-\sqrt{\rho\tau}\over 2n_{\rho,\tau}^2}\phe(n_{\rho,\tau}z+w)\cr
&\phantom{-\rho z^2+(1-\tau)zw+{\tau\over\rho}\biggl[\quad}\left.+
{1\over m_{\rho,\tau}^2}\psi(m_{\rho,\tau}z+w)-{1\over n_{\rho,\tau}^2}\psi(n_{\rho,\tau}z+w)\right],
\cr
&\quad(1-\rho)zw-\tau w^2+{\sqrt{\rho+\tau-1}+\sqrt{\rho\tau}\over 2}\phe(m_{\rho,\tau}z+w)\cr
&\phantom{\quad(1-\rho)zw-\tau w^2\,}
+{\sqrt{\rho+\tau-1}-\sqrt{\rho\tau}\over 2}\phe(n_{\rho,\tau}z+w)\cr
&\phantom{\quad(1-\rho)zw-\tau w^2\,}+\psi(m_{\rho,\tau}z+w)-\psi(n_{\rho,\tau}z+w)
\Biggr)\;,\cr}
$$ 
where the square roots of $\rho\tau$ and of $\rho+\tau-1$ are chosen so that their product is
equal to the previously chosen square root of $\rho\tau(\rho+\tau-1)$, and $\phe$,~$\psi\in\C[\![\zeta]\!]$ are power series of order at least~3.

\smallsect 3. Germs tangent to the identity

In this section we shall assume $n=\mu=2$ and $\Lambda=I$, that is we shall be interested
in 2-dimensional germs tangent to the identity of order~2. We shall keep using the notations introduced 
in the previous section. It should be recall that in his monumental work [\'E1] (see [\'E2] for
a survey) \'Ecalle studied the formal classification of germs tangent to the identity in dimension~$n$, giving a complete set of formal invariants for germs satisfying a generic condition:
the existence of at least one non-degenerate characteristic direction (an eigenradius, in
\'Ecalle's terminology). A {\sl characteristic direction} of a germ tangent to the identity $F$ is
a non-zero direction $v$ such that $F_\mu(v)=\lambda v$ for some $\lambda\in\C$, where
$F_\mu$ is the first (nonlinear) non-vanishing term in the homogeneous expansion of~$F$. The characteristic
direction $v$ is {\sl degenerate} if $\lambda=0$.

For this reason, we decided to discuss here the cases {\it without} non-degenerate
characteristic directions, that is the cases $(1_{00})$, $(1_{10})$ and $(2_{001})$, that cannot be
dealt with \'Ecalle's methods. Furthermore, we shall also study the somewhat special case
$(\infty)$, where all directions are characteristic; and we shall examine in detail case $(2_{10\rho})$, 
where interesting second-order resonance phenomena appear.

When $\Lambda=I$ the operator $L=L_{F_2,\Lambda}$ is given by
$$
L(H)=\Jac(H)\cdot F_2-\Jac(F_2)\cdot H\;.
$$
In particular, $L(F_2)=O$ always; therefore we cannot apply Proposition~\dinf\ (nor other 
similar conditions stated in [WZP2]), and we shall compute the second order normal form only.


\medbreak
\noindent$\bullet$ {\it Case $(\infty)$.}

\noindent In this case we have
$$
L(u_{d,j})=(d-2)u_{d+1,j+1}-v_{d+1,j}\qquad\hbox{and}\qquad L(v_{d,j})=(d-1)v_{d+1,j+1}
$$
for all $d\ge 2$ and $j=0,\ldots,d$. Therefore
$$
\Im L|_{\ca H^d}=
\cases{\Span{u_{d+1,2},\ldots,u_{d+1,d+1}, (d-2)u_{d+1,1}-v_{d+1,0},v_{d+1,1},\ldots,v_{d+1,d+1}}&for $d>2$,\cr
\Span{v_{3,0},\ldots,v_{3,3}}&for $d=2$.\cr}
$$
Thus
$$
(\Im L|_{\ca H^d})^\perp=
\cases{\Span{u_{d+1,0},(d+1)u_{d+1,1}+(d-2)v_{d+1,0}}&for $d>2$,\cr
\Span{u_{3,0},\ldots,u_{3,3}}&for $d=2$.\cr}
$$
It then follows that every formal power series of the form
$$
F(z,w)=(z+z^2+O_3,w+zw+O_3)
$$
has as second order normal form
$$
G(z,w)=\left(z+z^2+a_0z^3+a_1z^2w+a_2 zw^2+ \phe(w)+z\psi'(w),zw+w\psi'(w)-3\psi(w)\right)
$$
where $\phe\in\C[\![\zeta]\!]$ is a power series of order at least~3, $\psi\in\C[\![\zeta]\!]$ is a power series of order at least~4 and $a_0$,~$a_1$,~$a_2\in\C$.

\medbreak
\noindent$\bullet$ {\it Case $(1_{00})$.}

\noindent In this case we have
$$
L(u_{d,j})=(j-d)u_{d+1,j+2}+2v_{d+1,j+1}\qquad\hbox{and}\qquad L(v_{d,j})=(j-d)v_{d+1,j+2}
$$
for all $d\ge 2$ and $j=0,\ldots,d$. Therefore
$$
\Im L|_{\ca H^d}=\Span{2v_{d+1,1}-du_{d+1,2}, u_{d+1,3},\ldots,u_{d+1,d+1},v_{d+1,2},\ldots,v_{d+1,d+1}}\;,
$$
and thus
$$
(\Im L|_{\ca H^d})^\perp=\Span{u_{d+1,0},u_{d+1,1},v_{d+1,0},u_{d+1,2}+v_{d+1,1}}\;.
$$
It then follows that every formal power series of the form
$$
F(z,w)=(z+O_3,w-z^2+O_3)
$$
has as second order normal form
$$
G(z,w)=\left(z+w\phe_1(w)+z\phe_2(w)+z^2\psi(w),w-z^2+w\phe_3(w)+zw\psi(w)\right)\;,
$$
where $\phe_1$, $\phe_2$, $\phe_3\in\C[\![\zeta]\!]$ are power series of order at least~2, and $\psi\in\C[\![\zeta]\!]$ is a power series of
order at least~1.

\medbreak
\noindent$\bullet$ {\it Case $(1_{10})$.}

\noindent In this case we have
$$
L(u_{d,j})=(2-d)u_{d+1,j+1}-(d-j)u_{d+1,j+2}+2v_{d+1,j+1}+v_{d+1,j}
$$
and
$$
L(v_{d,j})=(1-d)v_{d+1,j+1}-(d-j)v_{d+1,j+2}
$$
for all $d\ge 2$ and $j=0,\ldots,d$. Therefore
$$
\Im L|_{\ca H^d}=
\cases{\Span{(2-d)u_{d+1,1}+v_{d+1,0},u_{d+1,2},\ldots,u_{d+1,d+1},v_{d+1,1},\ldots,v_{d+1,d+1}}&for $d>2$,\cr
\Span{v_{3,0}-2u_{3,2},u_{3,3},v_{3,1},v_{3,2},v_{3,3}}&for $d=2$,\cr}
$$
and thus
$$
(\Im L|_{\ca H^d})^\perp=
\cases{\Span{u_{d+1,0},(d+1)u_{d+1,1}+(d-2)v_{d+1,0}}&for $d>2$,\cr
\Span{u_{3,0},u_{3,1},3u_{3,2} + 2v_{3,0}}&for $d=2$.\cr}
$$
It then follows that every formal power series of the form
$$
F(z,w)=(z-z^2+O_3,w-z^2-zw+O_3)
$$
has as second order normal form
$$
{
G(z,w)=\bigl(z-z^2+\phe(w)+a_1zw^2+3a_2z^2w+z\psi'(w), w-z^2-zw+2a_2w^3+w\psi'(w)-3\psi(w)\bigr)\;,}
$$
where $\phe\in\C[\![\zeta]\!]$ is a power series of order at least~3, $\psi\in\C[\![\zeta]\!]$ is a power series of order at least~4, and $a_1$,~$a_2\in\C$.

%

\medbreak
\noindent$\bullet$ {\it Case $(2_{001})$.}

\noindent In this case we have
$$
L(u_{d,j})=(d-j)u_{d+1,j+1}-v_{d+1,j}\qquad\hbox{and}\qquad L(v_{d,j})=(d-j-1)v_{d+1,j+1}
$$
for all $d\ge 2$ and $j=0,\ldots,d$. It follows that 
$$
\Im L|_{\ca H^d}=\Span{du_{d+1,1}-v_{d+1,0},u_{d+1,2},\ldots,u_{d+1,d},v_{d+1,1},\ldots,v_{d+1,d+1}}
$$
and hence
$$
(\Im L|_{\ca H^d})^\perp=\Span{u_{d+1,0},u_{d+1,d+1},(d+1)u_{d+1,1}+dv_{d+1,0}}\;.
$$
It then follows that every formal germ of the form
$$
F(z,w)=(z+O_3,w+zw+O_3)
$$
has as second order normal form
$$
G(z,w)=\bigl(z+\phe_1(z)+\phe_2(w)+z\psi'(w),zw+w\psi'(w)-\psi(w)\bigr)
$$
where $\phe_1$, $\phe_2$, $\psi\in\C[\![\zeta]\!]$ are power series of order at least~3.

%

\medbreak
\noindent$\bullet$ {\it Case $(2_{10\rho})$.}

\noindent In this case we have
$$
L(u_{d,j})=(d-j-d\rho+2\rho) u_{d+1,j+1}+(\rho-1)v_{d+1,j}\quad\hbox{and}\quad L(v_{d,j})=(d-j-d\rho+\rho-1)v_{d+1,j+1}
\neweq\eqquno
$$
for all $d\ge 2$ and $j=0,\ldots,d$. Here we shall see the resonance phenomena we 
mentioned at the beginning of this section: for some values of~$\rho$ the dimension of the kernel of $L|_{\ca H^d}$ increases, and in some cases
we shall end up with a normal form depending on power series evaluated in monomials
of the form $z^{b-a}w^a$.

Let us put
$$
E_d=\left\{{d-j-1\over d-1}\biggm| j=0,\ldots, d\right\}\setminus\{0\}\quad\hbox{and}\quad
F_d=\left\{{d-j\over d-2}\biggm| j=0,\ldots, d-1\right\}
$$
(we are excluding 0 because $\rho\ne0$ by assumption),
where $E_d$ is defined for all $d\ge 2$ whereas $F_d$ is defined for all $d\ge 3$, and set
$$
\ca E=\bigcup_{d\ge 2} E_d=\bigl((0,1]\cap\Q\bigr)\cup\left\{-{1\over n}\biggm| n\in\N^*\right\}
$$
and
$$
\ca F=\bigcup_{d\ge 3} F_d=\bigl((0,1]\cap\Q\bigr)\cup\left\{1+{1\over n},1+{2\over n}\biggm| n\in\N^*\right\}.
$$
So $\ca E$ is the set of $\rho \in \C^*$ such that $L(v_{d, j})=0$ for some $d\ge 2$ and $0\le j\le d$, while $\ca F$ is the set of $\rho \in \C^*$ such that $L(u_{d, j})=(\rho -1) v_{d+1, j}$ for some $d\ge 3$ and $0\le j\le d-1$. 

\sm
Let us first discuss the non-resonant case, when $\rho\not\in\ca E\cup\ca F$. Then none of the
coefficients in \eqquno\ vanishes, and thus
$$
\Im L|_{\ca H^2}=\Span{2 u_{3,1}+(\rho-1)v_{3,0},u_{3,2},v_{3,1},v_{3,2},v_{3,3}}
$$
and
$$
\Im L|_{\ca H^d}=\Span{(d-d\rho+2\rho) u_{d+1,1} + (\rho -1) v_{d+1,0}, u_{d+1, 2},\ldots,u_{d+1,d+1},v_{d+1,1},\ldots,v_{d+1,d+1}}\;,
$$
for $d\ge 3$, and hence
$$
(\Im L|_{\ca H^d})^\perp=\cases{
\Span{u_{d+1,0},(1-\rho)(d+1)u_{d+1,1} + (d(1-\rho) + 2\rho) v_{d+1,0}}&for $d\ge 3$,\cr
\Span{u_{3,0},u_{3,3}, 3(1-\rho)u_{3,1} + 2v_{3,0}}& for $d=2$.\cr}
$$
It then follows that every formal germ of the form
$$
F(z,w)=(z-\rho z^2+O_3,w+(1-\rho)zw+O_3)
$$
with $\rho\not\in\ca E\cup\ca F$ (and $\rho\ne 0$) has as second order normal form
$$
G(z,w)=\bigl(z-\rho z^2+az^3+\phe(w) + (1-\rho) z\psi'(w), w + (1-\rho)zw + (1-\rho) w\psi'(w) + (3\rho -1) \psi(z)\bigr)\;,
$$
where $\phe$,~$\psi\in\C[\![\zeta]\!]$ are power series of order at least~3, and $a\in\C$.

Assume now $\rho\in\ca F\setminus\ca E$. Then $L(v_{d,j})\ne O$ always, and thus
$v_{d+1,j}\in\Im L|_{\ca H^d}$ for all $d\ge 2$ and all $j=1,\ldots, d+1$. Since $\rho>1$, if $d>2$ it also follows 
that $u_{d+1,j+1}\in\Im L|_{\ca H^d}$ for $j=1,\ldots,d$. 

Now, if $\rho=1+(1/n)$ then
$$
{d\over d-2}=\rho\quad\Longleftrightarrow\quad d=2(n+1)\;,
$$
and
$$
{d-1\over d-2}=\rho\quad\Longleftrightarrow\quad d=n+2\;.
$$
Taking care of the case $d=2$ separately, we then have
$$
\displaylines{
\Im L|_{\ca H^d}\hfill\cr
=\cases{
\Span{(d-d\rho+2\rho) u_{d+1,1} + (\rho -1) v_{d+1,0}, u_{d+1, 2},\ldots,u_{d+1,d+1},v_{d+1,1},\ldots,v_{d+1,d+1}}\kern-55pt&\cr
&\kern-44pt for $d\ge 3$, $d\ne n+2$, $2(n+1)$,\cr
\Span{u_{d+1,1}+(\rho-1) v_{d+1,0}, u_{d+1,3},\ldots,u_{d+1,d+1},v_{d+1,1},\ldots,v_{d+1,d+1}}&for $d= n+2$,\cr
\Span{u_{d+1,2},\ldots,u_{d+1,d+1},v_{d+1,0},\ldots,v_{d+1,d+1}}&for $d= 2(n+1)$,\cr
\Span{2u_{3,1}+ (\rho -1) v_{3,0},u_{3,2},v_{3,1},v_{3,2},v_{3,3}}&for $d=2$,
\cr}\hfill
\cr}
$$
and hence
$$
\displaylines{
(\Im L|_{\ca H^d}\!)^\perp\hfill\cr
=\cases{
\Span{u_{d+1,0},(1-\rho)(d+1) u_{d+1, 1} + (d(1-\rho)+2\rho)v_{d+1,0}}&for $d\ge 3$, $d\ne n+2$, $2(n+1)$,\cr
\Span{u_{d+1,0},u_{d+1,2},(1-\rho)(d+1)u_{d+1, 1} + v_{d+1, 0}}&for $d= n+2$,\cr
\Span{u_{d+1,0},u_{d+1,1}}&for $d= 2(n+1)$,\cr
\Span{u_{3,0},u_{3,3}, 3(1-\rho)u_{3,1} + 2 v_{3,0}}&for $d=2$.
\cr}\hfill\cr}
$$
It then follows that every formal germ of the form
$$
F(z,w)=\left(z-\left(1+{1\over n}\right) z^2+O_3,w-{1\over n}zw+O_3\right)
$$
with $n\in\N^*$ has as second order normal form
$$
\eqalign{
G(z,w)=
&\left(z-\left(1+{1\over n}\right) z^2 + \phe(w) + (1-\rho) z\psi'(w) + a_0 z^3+a_1 z^2 w^{n+1},\right.\cr 
&\quad\left.w-{1\over n}zw + (1-\rho) w\psi'(w) + (3\rho -1) \psi(w)\right)\;,
}
$$
where $\phe$,~$\psi\in\C[\![\zeta]\!]$ are power series of order at least~3, and $a_0$, $a_1\in\C$.

If instead $\rho=1+(2/m)$ with $m$ odd (if $m$ is even we are again in the previous case) then 
$$
{d\over d-2}=\rho\quad\Longleftrightarrow\quad d=m+2\;,
$$
whereas ${d-1\over d-2}\ne\rho$ always. Hence
$$
\displaylines{
\Im L|_{\ca H^d}\hfill\cr
=\cases{
\Span{(d-d\rho+2\rho) u_{d+1,1} + (\rho -1) v_{d+1,0},u_{d+1,2},\ldots,u_{d+1,d+1},v_{d+1,1},\ldots,v_{d+1,d+1}}\kern-55pt&\cr
&\kern-30pt for $d\ge 3$, $d\ne m+2$,\cr
\Span{u_{d+1,2},\ldots,u_{d+1,d+1},v_{d+1,0},\ldots,v_{d+1,d+1}}&for $d= m+2$,\cr
\Span{2u_{3,1}+(\rho -1) v_{3,0},u_{3,2},v_{3,1},v_{3,2},v_{3,3}}&for $d=2$,
\cr}\hfill
\cr}
$$
and thus
$$
(\Im L|_{\ca H^d})^\perp=\cases{
\Span{u_{d+1,0},(1-\rho) (d+1)u_{d+1,1} + (d-d\rho+2\rho)  v_{d+1,0}}&for $d\ge 3$, $d\ne m+2$,\cr
\Span{u_{d+1,0},u_{d+1,1}}&for $d= m+2$,\cr
\Span{u_{3,0},u_{3,3},3(1-\rho) u_{3,1} + 2v_{3,0}}&for $d=2$.
\cr}
$$
It then follows that every formal germ of the form
$$
F(z,w)=\left(z-\left(1+{2\over m}\right) z^2+O_3,w-{2\over m}zw+O_3\right)
$$
with $m\in\N^*$ odd has as second order normal form
$$
\eqalign{
G(z,w)=
&\left(z-\left(1+{2\over m}\right) z^2+\phe(w)+a_0 z^3+ (1-\rho)z(w\psi'(w) + \psi(w)),\right.\cr
&\quad\left.w-{2\over m}zw + (1-\rho) w^2 \psi'(w) + 2\rho w\psi(w)\right)\;,
}
$$
where $\phe\in\C[\![\zeta]\!]$ is a power series of order at least~3, $\psi\in\C[\![\zeta]\!]$ is a power series of order at least~2, and $a_0\in\C$.

Now let us consider the case $\rho=-1/n\in\ca E\setminus\ca F$. In this case the coefficients
in the expression of $L(u_{d,j})$ are always different from zero (with the exception of $d=j=2$),
whereas 
$$
d-j-d\rho+\rho-1=0\quad\Longleftrightarrow\quad j=d=n+1\;.
$$
It follows that
$$
\displaylines{
\Im L|_{\ca H^d}\hfill\cr
=\cases{
\Span{(d-d\rho+2\rho) u_{d+1,1} + (\rho -1) v_{d+1,0},u_{d+1,2},\ldots,u_{d+1,d+1},v_{d+1,1},\ldots,v_{d+1,d+1}}\kern-57pt&\cr
&for $d\ge 3$, $d\ne n+1$,\cr
\Span{(d-d\rho+2\rho) u_{d+1,1} + (\rho -1) v_{d+1,0},u_{d+1,2},\ldots,u_{d+1,d+1},v_{d+1,1},\ldots,v_{d+1,d}}\kern-57pt&\cr
&for $d= n+1$,\cr
\Span{2u_{3,1}+(\rho-1)v_{3,0},u_{3,2},v_{3,1},v_{3,2},v_{3,3}}&for $d=2$,
\cr}\hfill\cr}
$$
and thus
$$
\displaylines{
(\Im L|_{\ca H^d}\!)^\perp\hfill\cr
=\cases{
\Span{u_{d+1,0},(1-\rho)(d+1)u_{d+1,1} + (d-d\rho+2\rho) v_{d+1,0}}&for $d\ge 3$, $d\ne n+1$,\cr
\Span{u_{d+1,0},v_{d+1,d+1},(1-\rho)(d+1)u_{d+1,1} + (d-d\rho+2\rho) v_{d+1,0}}&for $d= n+1$,\cr
\Span{u_{3,0},u_{3,3},3(1-\rho)u_{3,1} + 2v_{3,0}}&for $d=2$.
\cr}\hfill\cr}
$$
It then follows that every formal germ of the form
$$
F(z,w)=\left(z+{1\over n}z^2+O_3,w+\left(1+{1\over n}\right)zw+O_3\right)
$$
with $n\in\N^*$ has as second order normal form
$$
\eqalign{
G(z,w)&=\left(z+{1\over n} z^2+\phe(w)+a_0 z^3+(1-\rho)z(w\psi'(w) + \psi(w)),\right.\cr
&\left.\qquad w+\left(1+{1\over n}\right)zw+\psi(z)+a_1z^{n+2} + (1-\rho)w^2\psi'(w) +2\rho w \psi(w) \right)\;,\cr}
$$
where $\phe\in\C[\![\zeta]\!]$ is a power series of order at least~3, $\psi\in\C[\![\zeta]\!]$ is a power series of order at least~2, and $a_0$, $a_1\in\C$.

Let us now discuss the extreme case
$\rho=1$. It is clear that
$$
\Im L|_{\ca H^d}=\Span{u_{d+1,1},u_{d+1,2},u_{d+1,4},\ldots,u_{d+1,d+1},v_{d+1,2},\ldots,v_{d+1,d+1}}\;,
$$
and hence
$$
(\Im L|_{\ca H^d})^\perp=\Span{u_{d+1,0},u_{d+1,3},v_{d+1,0},v_{d+1,1}}\;,
$$
It then follows that every formal germ of the form
$$
F(z,w)=\left(z-z^2+O_3,w+O_3\right)
$$
has as second order normal form
$$
G(z,w)=\bigl(z- z^2+\phe_1(w)+z^3\psi(w),
w+\phe_2(w)+z\phe_3(w)\bigr)\;,
$$
where $\phe_1$,~$\phe_2\in\C[\![\zeta]\!]$ are power series of order at least~3, $\phe_3\in\C[\![\zeta]\!]$ is a power series of order at least~2, and $\psi\in\C[\![\zeta]\!]$ is a power series.

We are left with the case $\rho\in(0,1)\cap\Q$. Write $\rho=a/b$ with $a$,~$b\in\N$ coprime and $0<a< b$. Now
$$
d-j-1-{a\over b}(d-1)=0\quad\Longleftrightarrow\quad j={(d-1)(b-a)\over b}\;;
$$
since $a$ and $b$ are coprime, this happens if and only if $d=b\ell+1$ and $j=(b-a)\ell$ for some $\ell\ge 1$. Analogously,
$$
d-j-{a\over b}(d-2)=0\quad\Longleftrightarrow\quad j=d-{a(d-2)\over b}\;;
$$
again, being $a$ and $b$ coprime, this happens if and only if $d=b\ell+2$ and $j=(b-a)\ell+2$ for some $\ell\ge 0$.
It follows that
$$
\displaylines{
\Im L|_{\ca H^d}\hfill\cr
=\cases{
\Span{(d-d\rho+2\rho) u_{d+1,1} + (\rho -1) v_{d+1,0},u_{d+1,2},\ldots,u_{d+1,d+1},v_{d+1,1},\ldots,v_{d+1,d+1}}\kern-60pt&\cr
&\kern-55pt for $d\ge 3$, $d\not\equiv 1$, 2 (mod $b$)\cr
\hbox{\rm Span}\,((d-d\rho+2\rho) u_{d+1,1} + (\rho -1) v_{d+1,0},u_{d+1,2},\ldots,\widehat{u_{d+1,(b-a)\ell+2}},\ldots,u_{d+1,d+1},\kern-60pt&\cr
\kern50pt v_{d+1,1},\ldots,\widehat{v_{d+1,(b-a)\ell+1}},\ldots,v_{d+1,d+1},{a\over b}u_{d+1,(b-a)\ell+2}-\left({a\over b}-1\right)v_{d+1,(b-a)\ell + 1})\kern-60pt&\cr
&for $d=b\ell +1$,\cr
\hbox{\rm Span}\,((d-d\rho+2\rho) u_{d+1,1} + (\rho -1) v_{d+1,0},u_{d+1,2},\ldots,\widehat{u_{d+1,(b-a)\ell+3}},\ldots, u_{d+1,d+1},\kern-60pt&\cr
\kern50pt v_{d+1,1},\ldots,v_{d+1,d+1})&for $d=b\ell+ 2$,\cr
\Span{2u_{3,1}+(\rho-1)v_{3,0},u_{3,2},v_{3,1},v_{3,2},v_{3,3}}&for $d=2$,
\cr}\hfill\cr}
$$
(where the hat indicates that that term is missing from the list), and thus
$$
\displaylines{
(\Im L|_{\ca H^d})^\perp\hfill\cr
=\cases{
\Span{u_{d+1,0},(1-\rho) u_{d+1, 1} + (d-d\rho+2\rho)v_{d+1,0}}&\kern-58pt for $d\ge 3$, $d\not\equiv 1$, 2 (mod $b$),\cr
\hbox{\rm Span}\,(u_{d+1,0},(1-\rho) u_{d+1, 1} + (d-d\rho+2\rho)v_{d+1,0},&\cr
\kern40pt (b-a)(a\ell+1)u_{d+1,(b-a)\ell+2}+a\bigl((b-a)\ell+2\bigr)v_{d+1,(b-a)\ell+1})&for $d=b\ell +1$,\cr
\Span{u_{d+1,0},u_{d+1,(b-a)\ell+3},(1-\rho) u_{d+1, 1} + (d-d\rho+2\rho)v_{d+1,0}}&for $d=b\ell+ 2$,\cr
\Span{u_{3,0},u_{3,3},3(1-\rho) u_{3,1} + 2 v_{3,0}}&for $d=2$.
\cr}\hfill\cr}
$$
It then follows that every formal germ of the form
$$
F(z,w)=\left(z-{a\over b}z^2+O_3,w+\left(1-{a\over b}\right)zw+O_3\right)
$$
with $a/b\in(0,1)\cap\Q$ and $a$, $b$ coprime, has as second order normal form
$$
\displaylines{
G(z,w)\hfill\cr
=\left(z-{a\over b} z^2+\phe(w)+z^3\phe_0(z^{b-a}w^a)+(b-a){\de\over\de w}\bigl(
z^2w\chi(z^{b-a}w^a)\bigr)+ \left(1-{a\over b}\right) z (w\psi'(w) + \psi(w) \bigr),\right.\cr
\left.\qquad w+\left(1-{a\over b}\right)zw+ a{\de\over\de z}\bigl(z^2w
\chi(z^{b-a}w^a)\bigr)+\left(1-{a\over b}\right) w^2\psi'(w) + 2{a\over b}w\psi(w) \right)\;,\cr}
$$
where $\phe$,~$\psi\in\C[\![\zeta]\!]$ are power series of order at least~3, and $\phe_0$,~$\chi\in\C[\![\zeta]\!]$ are power series of order at least~1.

%
%
%
%
%
%

%

\bigbreak
\noindent {\bf References.}\setref{AFGG}
\medskip

\art A M. Abate: Holomorphic classification of $2$-dimensional quadratic maps tangent to
the identity! S\=urikaisekiken\-ky\=usho K\=oky\=uroku! 1447 2005 1-14

\coll A2 M. Abate: Discrete holomorphic local dynamical systems!
Holomorphic dynamical systems! G.Gentili, J. Gu\'enot, G. Patrizio eds., Lect. Notes
in Math. {\bf 1998,} Springer, Berlin, 2010, pp. 1--55


\art AT1 M. Abate, F. Tovena: Formal classification of holomorphic maps tangent
to the identity! Disc. Cont. Dyn. Sys.! Suppl. 2005 1-10

\pre AT2 M. Abate, F. Tovena: Poincar\'e-Bendixson theorems for meromorphic connections and homogeneous vector fields! To appear in J. Diff. Eq. DOI: 10.1016/j.jde.2011.05.031, arXiv:0903.3485! 2011

\art AFGG A. Algaba, E. Freire, E. Gamero. C. Garcia: Quasi-homogeneous normal forms! J. Comput. Appl. Math.! 150 2003 193-216

\book Ar V.I. Arnold: Geometrical methods in the theory of ordinary differential equations! Springer Verlag, New York, 1988

\art B1 G.R. Belitskii: Invariant normal forms of formal series! Funct. Anal. Appl.! 13 1979 46-67

\art B2 G. R. Belitskii: Normal forms relative to a filtering action of a group! Trans. Moscow Math. Soc.! 1981, 1982 1-39

\art BS A. Baider, J. Sanders: Further reduction of the Takens-Bogdanov normal form! J. Diff. Eq.! 99 1992 205-244

%

\book C H. Cartan: Cours de calcul diff\'erentiel, Hermann, Paris! 1977

\art CD G.T. Chen, J.D. Dora: Normal forms for differentiable maps near a fixed point! Numer. Algorith.! 22 1999 213-230

%

\book \'E1 J. \'Ecalle: Les fonctions r\'esurgentes. Tome I\negthinspace
I\negthinspace I: L'\'equation du pont et la classification analytique
des objects locaux! Publ. Math. Orsay {\bf 85-05,} Universit\'e de
Paris-Sud, Orsay, 1985 

\coll \'E2 J. \'Ecalle: Iteration and analytic classification of local
diffeomorphisms of $\C^\nu$! Iteration theory and its
functional equations! Lect. Notes in Math. {\bf 1163,} Springer-Verlag,
Berlin, 1985, pp. 41--48

%

\art F E. Fischer: \"Uber die Differentiationsprozesse der Algebra! J. f\"ur Math.! 148 1917 1-78

\art G G. Gaeta: Further reduction of Poincar\'e-Dulac normal forms in symmetric systems! Cubo! 9 2007 1-11

\art KOW H. Kobuki, H. Oka, D. Wang: Linear grading function and further reduction
of normal forms! J. Diff. Eq.! 132 1996 293-318

\art LS E. Lombardi, L. Stolovitch: Normal forms of analytic perturbations of quasihomogeneous vector fields: Rigidity, analytic invariant sets and exponentially small approximation! Ann. Scient. Ec. Norm. Sup.! 43 2010 659-718

%

\book Mu1 J. Murdock: Normal forms and unfoldings for local dynamical systems! Springer Verlag,
Berlin, 2003

\art Mu2 J. Murdock: Hypernormal form theory: foundations and algorithms! J. Diff. Eq.! 205 2004 424-465

\art R1 J. Raissy: Torus actions in the normalization problem! J. Geom. Anal.! 20 2010 472-524

\book R2 J. Raissy: {\sl Brjuno conditions for linearization in presence of resonances}! ArXiv:0911.4341. 
To appear in  {\bf Asymptotics in Dynamics, Geometry and PDE's; Generalized Borel Summation}, Eds. O. Costin, F. Fauvet, F. Menous and D. Sauzin, Edizioni Della Normale, Pisa, 2010

\art R\"u H. R\"ussmann: Stability of elliptic fixed points of analytic area-preserving mappings 
under the Brjuno condition! Ergodic Theory Dynam. Systems! 22 2002 1551-1573

\art WZP1 D. Wang, M. Zheng, J. Peng: Further reduction of normal forms of formal maps! C. R. Math. Acad. Sci. Pari! 343 2006 657-660

\art WZP2 D. Wang, M. Zheng, J. Peng: Further reduction of normal forms and unique normal forms of smooth maps! Internat. J. Bifur. Chaos! 18 2008 803-825

%

\bye